\RequirePackage{fix-cm}
\documentclass[smallcondensed]{svjour3}     % onecolumn (ditto)
\smartqed  % flush right qed marks, e.g. at end of proof
\usepackage{graphicx}
\usepackage{latexsym,amsfonts,amsbsy,amssymb}
\usepackage{amsmath,amscd}
\usepackage{hyperref}
\usepackage[linesnumbered, ruled]{algorithm2e}
\usepackage{float}
\usepackage{cite}
\usepackage{bm}
\usepackage{upgreek}
\usepackage{url}
\usepackage[page]{appendix}
\journalname{Journal of Scientific Computing}
\begin{document}

\title{On the Asymptotic Linear Convergence Speed of Anderson Acceleration Applied to ADMM%\thanks{Grants or other notes
%about the article that should go on the front page should be
%placed here. General acknowledgments should be placed at the end of the article.}
}
%\subtitle{Do you have a subtitle?\\ If so, write it here}

\titlerunning{On the Asymptotic Linear Convergence Speed of AA-ADMM}        % if too long for running head

\author{Dawei Wang         \and
        Yunhui He          \and
        Hans De Sterck
}

%\authorrunning{Dawei Wang \and
%               Yunhui He  \and
%               Hans De Sterck} % if too long for running head

\institute{D. Wang\at
              Department of Applied Mathematics, University of Waterloo, 200 University Ave W, ON. N2L3G1, Canada\\
              \email{dawei.wang@uwaterloo.ca}           %  \\
%             \emph{Present address:} of F. Author  %  if needed
           \and
           Y. He \at
              Department of Applied Mathematics, University of Waterloo, 200 University Ave W, ON. N2L3G1, Canada\\
              \email{yunhui.he@uwaterloo.ca}           %  \\
%             \emph{Present address:} of F. Author  %  if needed
           \and
           H. De Sterck (Corresponding author)\at
              Department of Applied Mathematics, University of Waterloo, 200 University Ave W, ON. N2L3G1, Canada\\
              \email{hdesterck@uwaterloo.ca}           %  \\
%             \emph{Present address:} of F. Author  %  if needed
}

\date{Received: date / Accepted: date}
% The correct dates will be entered by the editor

\maketitle

\begin{abstract}
Empirical results show that Anderson acceleration (AA) can be a powerful
mechanism to improve the asymptotic linear convergence speed of the Alternating
Direction Method of Multipliers (ADMM) when ADMM by itself converges linearly.
However, theoretical results to quantify this improvement do not exist yet.
In this paper we explain and quantify this improvement in linear asymptotic convergence
speed for the special case of a stationary version of AA applied to ADMM. We do so by considering
the spectral properties of the Jacobians of ADMM and the stationary version of AA evaluated
at the fixed point, where the coefficients of the stationary AA method are computed such that
its asymptotic linear convergence factor is optimal.
The optimal linear convergence factors of this stationary AA-ADMM method are
computed analytically or by optimization, based on previous work on optimal stationary AA
acceleration. Using this spectral picture and those analytical results, our approach provides
new insight into how and by how much the stationary AA method can improve the
asymptotic linear convergence factor of ADMM.
Numerical results also indicate that the optimal linear convergence factor
of the stationary AA methods gives a useful estimate for the asymptotic linear
convergence speed of the non-stationary AA method that is used in practice.
%We explain how Anderson Acceleration (AA) speeds up the Alternating Direction
%Method of Multipliers (ADMM), for the case where ADMM by itself converges linearly.
%We do so by considering the spectral properties of the Jacobians
%of ADMM and a stationary version of AA evaluated at the fixed point,
%where the coefficients of the stationary version are computed such that
%its asymptotic linear convergence factor is optimal.
%Numerical tests show that this allows us to quantify the improved
%linear asymptotic convergence speed of AA-ADMM as compared to the
%convergence factor of ADMM used by itself.
%This way of estimating AA-ADMM convergence speed is useful because
%there are no known convergence bounds for AA with finite window size
%that would allow quantification of this improvement in asymptotic convergence
%speed.
\keywords{Anderson acceleration \and ADMM \and asymptotic linear convergence speed \and machine learning}
% \PACS{PACS code1 \and PACS code2 \and more}
\subclass{65K10}
\end{abstract}

%=============================================================
\section{Introduction}
\label{intro}
%=============================================================
In this paper, we consider the constrained optimization problem
%-----------------------------------------
\begin{equation}\label{eqn:opt-model}
\begin{aligned}
    \min_{\mathbf x,\mathbf z}\; &f(\mathbf x,\mathbf z)=f_1(\mathbf x) + f_2(\mathbf z), \\
    \text{s.t. }\; &\mathbf A\mathbf x + \mathbf B\mathbf z = \mathbf b,
\end{aligned}
\end{equation}
%-----------------------------------------
where $\mathbf x\in\mathbb R^{n_1}$, $\mathbf z\in\mathbb R^{n_2}$ are optimization variables, $\mathbf b\in\mathbb R^{n_b}$ is a known vector of data, $f_1 : \mathbb R^{n_1}\rightarrow\mathbb R$, $f_2 : \mathbb R^{n_2}\rightarrow\mathbb R$ are the objective functions, and $\mathbf A\in\mathbb R^{n_b\times n_1}, \, \mathbf B\in\mathbb R^{n_b\times n_2}$ are linear operators. Many optimization problems in data science and machine learning can be cast into this form.

We consider the well-known Alternating Direction Method of Multipliers (ADMM) \cite{boyd2011distributed}
for solving problem (\ref{eqn:opt-model}), and we apply Anderson acceleration (AA) \cite{anderson1965iterative}
to accelerate the convergence of ADMM. In particular, we consider problems where ADMM by itself would converge
linearly with a linear asymptotic convergence factor $\rho_{ADMM}$, and we are interested in explaining and
quantifying how and by how much the combined AA-ADMM method would improve the asymptotic convergence compared
to $\rho_{ADMM}$. In recent papers it has indeed been observed numerically that AA may speed up the
convergence of ADMM and related methods substantially \cite{zhang2019accelerating,fu2019anderson,mai2019anderson},
but there are no known convergence bounds for AA with finite window size that would allow quantification of this improvement in linear asymptotic convergence speed. 

Since the analysis of convergence acceleration by AA with finite window size has so far proven
intractable, we investigate in this paper the simplified case of convergence acceleration of ADMM
by  a stationary version of AA (sAA), where the sAA coefficients are determined in a way that
optimizes the asymptotic linear convergence factor of the stationary sAA-ADMM method,
given the spectral properties of the Jacobian of the ADMM update at the fixed point.
We will demonstrate how the spectral properties
of the ADMM and optimal sAA-ADMM Jacobians can be used to explain how and by how much
the sAA nonlinear convergence acceleration method can accelerate the asymptotic convergence
of ADMM.
We use the theoretical results that were introduced in \cite{desterck2020}
for analyzing convergence acceleration by stationary versions of
AA and the closely related nonlinear GMRES (NGMRES) method,
which were applied in \cite{desterck2020} to the acceleration of the
Alternating Least Squares (ALS)  method to compute canonical tensor decompositions.
AA (in its NGMRES form) was first applied to accelerate the convergence of ALS
for the nonconvex canonical tensor decomposition problem in 2012
\cite{sterck2012nonlinear}.
We use the theoretical results from \cite{desterck2020} on optimal sAA coefficients to
compute the optimal sAA-ADMM asymptotic convergence factor, $\rho^*_{sAA-ADMM}$.
We will also report on numerical tests indicating that the optimal stationary
$\rho^*_{sAA-ADMM}$ factors provide a useful estimate for the improved asymptotic linear convergence
speed of applying the non-stationary AA method that is used in practice to ADMM.

%+++++++++++++++++++++++++++++++++++++++++++++++++
\subsection{Alternating Direction Method of Multipliers}
%+++++++++++++++++++++++++++++++++++++++++++++++++
Extensive research has shown that ADMM is an effective tool for solving (\ref{eqn:opt-model}), and can be competitive with the best known methods for some problems \cite{boyd2011distributed},
in particular also when accelerated by AA \cite{zhang2019accelerating,fu2019anderson,mai2019anderson}. To present ADMM for solving (\ref{eqn:opt-model}), we first need to define the augmented Lagrangian
\begin{equation}\label{eqn:aug_lagrangian}
    L_{\rho}(\mathbf x,\mathbf z,\mathbf y) = f_1(\mathbf x) + f_2(\mathbf z) + \mathbf y^T(\mathbf A \mathbf x  + \mathbf B \mathbf z- \mathbf  b) + \frac{\rho}{2}||\mathbf A \mathbf x  + \mathbf B \mathbf z- \mathbf  b||^2_2,
\end{equation}
where $\mathbf y$ is the Lagrange multiplier, and $\rho > 0$ is a penalty parameter. ADMM then solves the original problem by performing alternating minimization of the augmented Lagrangian with respect to variables $\mathbf x$ and $\mathbf z$ and computes the sub-problems
\begin{equation*}
\begin{cases}
    \mathbf x_{k+1} = \text{argmin}_{\mathbf x} L_{\rho}(\mathbf x, \mathbf z_k, \mathbf y_k),\\
    \mathbf z_{k+1} = \text{argmin}_{\mathbf z}L_{\rho}(\mathbf x_{k+1},\mathbf z, \mathbf y_k), \\
    \mathbf y_{k+1} = \mathbf y_k + \rho(\mathbf A\mathbf x_{k+1}+\mathbf B \mathbf z_{k+1}-\mathbf b),
\end{cases}
\end{equation*}
given initial approximations $\mathbf z_0$ and $\mathbf y_0$. It is often more convenient to write the augmented Lagrangian (\ref{eqn:aug_lagrangian}) in an equivalent scaled form by replacing $\frac{1}{\rho}\mathbf y$ with $\mathbf u$
\begin{equation}\label{eqn:scaled_aug_lagrangian}
    L_{\rho}(\mathbf x,\mathbf z,\mathbf u) = f_1(\mathbf x) + f_2(\mathbf z) + \frac{\rho}{2}||\mathbf A \mathbf x  + \mathbf B \mathbf z- \mathbf  b+\mathbf u||^2_2 - \frac{\rho}{2}||\mathbf u||^2_2.
\end{equation}
% observing that
% \begin{align*}
%     &\mathbf y^T(\mathbf A \mathbf x  + \mathbf B \mathbf z- \mathbf  b) + \frac{\rho}{2}||\mathbf A \mathbf x  + \mathbf B \mathbf z- \mathbf  b||^2_2\\
%     =&\frac{\rho}{2}||\mathbf A \mathbf x  + \mathbf B \mathbf z- \mathbf  b+\frac{1}{\rho}y||^2_2 - \frac{1}{2\rho}||y||^2_2\\
%     =&\frac{\rho}{2}||\mathbf A \mathbf x  + \mathbf B \mathbf z- \mathbf  b+\frac{1}{\rho}y||^2_2 - \frac{\rho}{2}||\frac{1}{\rho}y||^2_2.
% \end{align*}
Then the ADMM steps become
\begin{equation}
\begin{cases}
    \mathbf x_{k+1} = \text{argmin}_{\mathbf x} f_1(\mathbf x) + \frac{\rho}{2}||\mathbf A \mathbf x+\mathbf B \mathbf z_k-\mathbf b+\mathbf u_k||^2_2,\\
    \mathbf z_{k+1} = \text{argmin}_{\mathbf z}f_2(\mathbf z) + \frac{\rho}{2}||\mathbf A\mathbf x_{k+1}+\mathbf B \mathbf z-\mathbf b+\mathbf u_k||^2_2, \\
    \mathbf u_{k+1} = \mathbf u_k + \mathbf A\mathbf x_{k+1}+\mathbf B \mathbf z_{k+1}-\mathbf b,
\end{cases}
\label{eq:myADMM}
\end{equation}
given initial approximations $\mathbf z_0$ and $\mathbf u_0$.

The optimality conditions for problem (\ref{eqn:opt-model}) using ADMM are the primal feasibility
\begin{align}
    \mathbf A \mathbf x^* + \mathbf B \mathbf z^*-\mathbf b = 0,
\end{align}
and dual feasibility
\begin{align}
    0 \in \partial f_1(\mathbf x^*) + \mathbf A^T\mathbf y^*,\label{eqn:dual_feasibility_x}\\
    0 \in \partial f_2(\mathbf z^*) + \mathbf B^T\mathbf y^*,\label{eqn:dual_feasibility_z}
\end{align}
where $\mathbf x^*, \mathbf z^*, \mathbf y^*$ are the optimal solutions. It turns out that $\mathbf z_{k+1}$ and $\mathbf y_{k+1}$ always satisfy dual feasibility (\ref{eqn:dual_feasibility_z}), and the optimization step
for $\mathbf x_{k+1}$ implies \cite{boyd2011distributed}
\begin{equation*}
    \rho \mathbf A^T\mathbf B(\mathbf z_{k+1}-\mathbf z_k)\in \partial f_1(\mathbf x_{k+1}) + A^T\mathbf y_{k+1}.
\end{equation*}
This means that 
\begin{equation*}
    \mathbf r^p_{k+1} := \mathbf A\mathbf x_{k+1}+\mathbf B \mathbf z_{k+1}-\mathbf b
\end{equation*}
can be used as the primal residual at iteration $k+1$, and
\begin{equation*}
    \mathbf r^d_{k+1} := \rho \mathbf A^T\mathbf B(\mathbf z_{k+1}-\mathbf z_k)
\end{equation*}
can be used as the dual residual at iteration $k+1$. These two residuals converge to zero as ADMM proceeds \cite{boyd2011distributed}. 
%We will follow \cite{goldstein2014fast} and refer to
%\begin{equation}\label{eqn:combined_r}
%    r^{c}_{k+1} = \rho||A\mathbf x_{k+1}+\mathbf B \mathbf z_{k+1}-\mathbf b||^2 + \rho||B(\mathbf z_{k+1}-z_k)||^2
%\end{equation}
% this is actually a SQUARED combined residual! don't use this terminology
%as the combined residual. It is also proven in \cite{goldstein2014fast} that the ADMM iterations decrease the combined residual (\ref{eqn:combined_r}) monotonically. This is useful in safeguarding schemes that may be used for accelerated ADMM.

Although there are abundant results on the application of ADMM, studies on ADMM convergence rates are few until recently. When the objective functions $f_1$ and $f_2$ are convex (not requiring strong convexity, and possibly nonsmooth), the work in \cite{he20121,he2015non,davis2017faster} has shown an $\mathcal O(1/k)$ convergence rate under some additional assumptions. 
The work in \cite{hong2017linear, boley2012linear, nishihara2015general, lions1979splitting, deng2016global, davis2017faster} shows linear convergence of ADMM under strong convexity and rank conditions. More specifically, results in \cite{lions1979splitting} show that when $f$ is strongly convex and the composite constraint matrix $[A\; B]$ is row independent, then ADMM converges linearly to the unique minimizer. More recent work in \cite{boley2012linear, deng2016global} shows that when at least one of the component functions is strongly convex and has a Lipschitz-continuous gradient, and under certain rank conditions on the constraint matrices, some linear convergence results can be obtained for a subset of primal and dual variables in the ADMM algorithm. The often slow convergence of ADMM is one of the reasons that ADMM was not well-known until recently when large-scale distributed optimization became necessary.

%+++++++++++++++++++++++++++++++++++++++++++++++++
\subsection{Acceleration methods for ADMM}
%+++++++++++++++++++++++++++++++++++++++++++++++++
Results on accelerated versions of ADMM are even fewer. The most widely used acceleration technique is simple overrelaxation, which reliably reduces the total iteration count by a small factor \cite{ghadimi2014optimal}. A GMRES-accelerated ADMM is discussed in \cite{zhang2018gmres} for a quadratic objective, for which the ADMM iteration is linear.
In some sense, our paper is a nonlinear extension of the approach in \cite{zhang2018gmres} since AA is a nonlinear generalization of GMRES \cite{walker2011anderson,desterck2020}:
we consider nonlinear convergence acceleration by AA of general nonlinear ADMM iterations that converge linearly, and
\cite{zhang2018gmres} considers linear convergence acceleration by GMRES of specific linear
ADMM iterations.
For the case of Nesterov acceleration, which is a version of Anderson acceleration with window size one \cite{mitchell2020nesterov,desterck2020}, the only papers providing convergence rates for not necessarily differentiable convex functions are \cite{goldstein2014fast, kadkhodaie2015accelerated,francca2018dynamical,franca2018admm}, among which \cite{goldstein2014fast,kadkhodaie2015accelerated} show that under strong convexity assumptions Nesterov acceleration of ADMM has an optimal global convergence bound of $\mathcal O(1/k^2)$ in terms of the primal and dual residual norms. In \cite{francca2018dynamical} a dynamical system perspective was proposed for understanding ADMM and accelerated ADMM applied to the problem (\ref{eqn:opt-model}) with the constraint $\mathbf z = \mathbf A \mathbf x$. Using a nonsmooth Lyapunov analysis technique, they proved a convergence rate of $\mathcal O(1/k)$ for ADMM, and a convergence rate of $\mathcal O(1/k^2)$ for accelerated ADMM, under the assumption that $f_1$ and $f_2$ are both proper, lower semicontinuous and convex, and $A$ has full column rank. Following this work, more convergence rates of dynamical systems related to relaxed and accelerated variants of ADMM are given in \cite{franca2018admm}.

Work using Anderson acceleration (AA) applied to ADMM and related methods can be found in \cite{zhang2019accelerating,peng2018anderson,kadkhodaie2015accelerated,fu2019anderson,poon2019trajectory}, but no convergence rates are given that quantify convergence improvement. 
In this paper, we investigate acceleration of ADMM by the stationary version of AA (sAA)
that was first introduced in \cite{desterck2020} for the case that ADMM converges linearly, and we
determine optimal linear asymptotic convergence factors for the accelerated sAA-ADMM algorithm,
quantifying the convergence improvement relative to the linear asymptotic convergence factor of ADMM
used by itself.
We also provide numerical results indicating that these optimal sAA convergence factors
give a useful estimate of the asymptotic convergence improvement provided by the non-stationary
AA method that is used in practice.
 
%The spectral properties of ADMM and sAA-ADMM at the fixed point will explain
%how sAA accelerates the linear asymptotic convergence speed of ADMM,
%and we will see in numerical results that the linear asymptotic convergence speed of the
%optimal stationary AA methods gives a useful estimate for the asymptotic convergence
%improvement provided by the non-stationary AA method that is used in practice.

%-----------------------------------------------------------------
\subsubsection{Anderson acceleration for fixed-point iterations}
%-----------------------------------------------------------------
Consider fixed-point iteration (FPI)
\begin{equation}\label{eq:FPI}
    \mathbf x_{k+1} = \mathbf q(\mathbf x_k),
\end{equation}
where $\mathbf q: \mathbb R^{n}\rightarrow\mathbb R^{n}$ is the iteration function. The method of Anderson acceleration tries to improve convergence by taking
\begin{equation} \label{eqn:anderson_k-i_k-i-1}
    \mathbf x_{k+1} = \mathbf q(\mathbf x_k) + \sum_{i=0}^{m_k-1} \beta^{(k)}_i\left(\mathbf q(\mathbf x_{k-i})-\mathbf q(\mathbf x_{k-i-1})\right).
\end{equation}
where $m_k = \min\{m,k\}$ with some predefined window size $m \geq 0$, and the coefficients
$\beta^{(k)}_i$ are computed from optimization problem
\begin{equation}\label{eqn:anderson_coeff_k-i_k-i-1}
    \{\beta_i^{(k)}\} = \underset{\{\beta_i\}}{\text{argmin}}||\mathbf r(\mathbf x_k) + \sum_{i=0}^{m_k-1} \beta_i\left(\mathbf r(\mathbf x_{k-i})-\mathbf r(\mathbf x_{k-i-1})\right)||^2,
\end{equation}
where $\mathbf r(\mathbf x_k) =\mathbf x_k - \mathbf q(\mathbf x_k)$ is the residual of FPI (\ref{eq:FPI}) in iteration $k$. 
We refer to Anderson acceleration with window size $m$ as AA($m$).

It has been shown that Anderson acceleration is, in the linear case, essentially equivalent to the GMRES method for solving linear systems when $m=k$ \cite{walker2011anderson}. When $m = 0$, the un-accelerated FPI is recovered. 
The convergence of Anderson acceleration is not guaranteed. The work in \cite{toth2015convergence} shows that for linear problems, if the FPI is a contraction, global convergence can be proved. But for nonlinear problems, only local convergence can be shown under certain conditions. Global convergence properties can be improved by adding a safeguarding step to the algorithm
\cite{desterck2020,zhang2019accelerating,fu2019anderson,mitchell2020nesterov}.
However, we do not need a safeguarding step for the numerical tests with linear asymptotic convergence that we consider in this paper.

In \cite{desterck2020}, a stationary variant of AA is considered, which we call sAA, and is given by
\begin{equation} \label{eqn:anderson_k-i_k-i-1-stat}
    \mathbf x_{k+1} = \mathbf q(\mathbf x_k) + \sum_{i=0}^{m_k-1} \beta_i\left(\mathbf q(\mathbf x_{k-i})-\mathbf q(\mathbf x_{k-i-1})\right),
\end{equation}
where the $\beta_i$ are fixed for all iterations. We refer to sAA with window size $m$ as sAA($m$).
In \cite{desterck2020}, the constant sAA coefficients $\beta_i$ in (\ref{eqn:anderson_k-i_k-i-1-stat}) are computed
such that the asymptotic linear convergence factor of the sAA method is optimal, given knowledge of $\mathbf q'(\mathbf x)$
evaluated in the fixed point $\mathbf x^*$ (see Section \ref{sec:opt} for details). 
We use this approach in this paper to quantify the optimal asymptotic convergence
speed of sAA-ADMM compared to $\rho_{ADMM}$, and the spectral properties of $\mathbf q'(\mathbf x^*)$ provide insight into how sAA effectively accelerates ADMM, as will be discussed in Section \ref{sec:results}.
%Numerical results in \cite{desterck2020} show that this optimal sAA
%asymptotic convergence factor can be used to estimate the asymptotic convergence speed for AA with
%finite window size.

%-----------------------------------------------------------------
\subsubsection{Anderson acceleration applied to ADMM (AA-ADMM)}
%-----------------------------------------------------------------
When we use AA to accelerate ADMM, we can treat one iterate of ADMM as a FPI, that is,
the ADMM iteration of (\ref{eq:myADMM})
%\begin{equation*}
%\begin{cases}
%    \mathbf x_{k+1} = \text{argmin}_{\mathbf x} f_1(\mathbf x) + \frac{\rho}{2}||\mathbf A \mathbf x+\mathbf B \mathbf z_k-\mathbf b+\mathbf u_k||^2_2,\\
%    \mathbf z_{k+1} = \text{argmin}_{\mathbf X}f_2(\mathbf z) + \frac{\rho}{2}||\mathbf A\mathbf x_{k+1}+\mathbf B \mathbf z-\mathbf b+\mathbf u_k||^2_2, \\
%    \mathbf u_{k+1} =\mathbf u_k + \mathbf A\mathbf x_{k+1}+\mathbf B \mathbf z_{k+1}-\mathbf b.
%\end{cases}
%\end{equation*}
can be seen as a FPI
\begin{equation}\label{eqn:fpi-qzu}
    (\mathbf z_{k+1}, \mathbf u_{k+1}) = \mathbf q(\mathbf z_k, \mathbf u_k),
\end{equation}
given initial approximations $\mathbf z_0,\;\mathbf u_0$. Notice that $\mathbf x_{k+1}$ is only dependent on $\mathbf z_k$ and $\mathbf u_k$ and can be recovered from them anytime during the iteration, thus it is included implicitly and can be eliminated when ADMM is seen as a FPI \cite{zhang2019accelerating}. Moreover, if $\mathbf B$ is a nonsingular square matrix, since
\begin{equation*}
    \nabla f_2(\mathbf z_{k+1}) + \rho \mathbf B^T(\mathbf A\mathbf x_{k+1}+\mathbf B \mathbf z_{k+1}-\mathbf b+\mathbf u_k) = 0,
\end{equation*}
from the step of the $\mathbf z_{k+1}$ update, we get
\begin{equation*}
    \mathbf u_k + A\mathbf x_{k+1} + \mathbf B \mathbf z_{k+1} - \mathbf b = -\frac{1}{\rho}\mathbf B^{-T}\nabla f_2(\mathbf z_{k+1}),
\end{equation*}
and thus
\begin{equation*}
    \mathbf u_{k+1} = -\frac{1}{\rho}\mathbf B^{-T}\nabla f_2(\mathbf z_{k+1}).
\end{equation*}
Then, we can further simplify ADMM as a FPI of variable $\mathbf z$ only \cite{zhang2019accelerating}, i.e.,
\begin{equation}\label{eqn:fpi-qz}
    \mathbf z_{k+1} = q(\mathbf z_k).
\end{equation}
The other two variables $\mathbf x_{k+1}$ and $\mathbf u_{k+1}$ can be recovered from $\mathbf z_k$.
% and thus can be eliminated when the ADMM is seen as a FPI. 
These simplifications are not necessary, but they help avoid computational overhead and simplify implementation.

The rest of this paper is structured as follows. In Section \ref{sec:opt} we discuss the detailed theoretical results on stationary AA from \cite{desterck2020} that will be used in this paper to analyze the 
convergence acceleration of ADMM by sAA in Section \ref{sec:results}. Section \ref{sec:results}
will also numerically compare acceleration of ADMM by stationary and non-stationary AA.
Conclusions are formulated in Section \ref{sec:conc}.

%=============================================================
\section{Optimal asymptotic convergence speed of stationary AA applied to ADMM}
\label{sec:opt}
As we mentioned earlier, there is a lack of mathematical understanding of the improved asymptotic
convergence speed of AA with finite window size applied to FPI (\ref{eq:FPI}). In this section, we discuss
the theory from \cite{desterck2020} that quantifies how the
stationary version of AA can optimally accelerate the asymptotic convergence of a linearly converging FPI. 
We summarize the results from \cite{desterck2020} with small extensions 
in a form that is convenient for the purposes of this paper.
This theory focuses on the analysis of sAA with window size $m = 1$, and it assumes that the fixed-point iteration operator $\mathbf q(\cdot)$ is differentiable at the fixed point $\mathbf x^*$, and that the FPI converges root-linearly with linear convergence factor $\rho$ that is the spectral radius of $\mathbf q'(\mathbf x^*)$.

We will apply this theory in this paper to quantify the improved asymptotic convergence speed of 
the stationary version of AA applied to ADMM, compared to $\rho_{ADMM}$, in the case that ADMM by itself converges linearly. We will make the assumption that the ADMM iteration operator $\mathbf q(\cdot)$ is differentiable at $\mathbf x^*$.
It is worth mentioning that for the analysis we pursue, we only need to assume the differentiability of $\mathbf q(\cdot)$ in a neighborhood of the solution, and $\mathbf q(\cdot)$ does not need to be smooth elsewhere.
In fact, the objective function $f(\mathbf x, \mathbf z)$ in (\ref{eqn:opt-model}) may not be differentiable at the solution, but this does not necessarily preclude the ADMM
iteration operator $\mathbf q(\mathbf x)$ from being differentiable at the solution.
We elaborate on this in Appendix \ref{appendix:l1-norm}, and this means that our
approach of analyzing sAA-ADMM convergence based on the spectral properties
of $\mathbf q'(\mathbf x^*)$ may be applied to both differentiable and non-differentiable
objectives $f(\mathbf x, \mathbf z)$ in (\ref{eqn:opt-model}), as long as $\mathbf q(\cdot)$ 
is differentiable in a neighborhood of the solution and the asymptotic convergence of ADMM
by itself is linear.

%+++++++++++++++++++++++++++++++++++++++++++++++++
%\subsection{Stationary Anderson acceleration}
%+++++++++++++++++++++++++++++++++++++++++++++++++
%The theoretical analysis of the improved convergence speed of AA remains a challenging topic; one of the reasons is that the AA coefficients $\beta_i^{(k)}$ in (\ref{eqn:anderson_k-i_k-i-1}) change at every iteration and thus are hard to analyze. A recent study by De Sterck and He \cite{desterck2020} investigates the asymptotic convergence factor of AA applied to the steepest descent method and Alternating Least Squares by fixing the coefficients $\beta_i$ in stationary variant (\ref{eqn:anderson_k-i_k-i-1-stat}) of AA (sAA), and choosing the coefficients to be optimal in terms of the sAA convergence factor. Numerical results in \cite{desterck2020} find that the asymptotic convergence speed of
%AA with finite window size, where the $\beta_i^{(k)}$ are determined in each iteration $k$ in a locally optimal way,
%is similar to the asymptotic convergence speed of the stationary sAA method with globally optimal coefficients.
%As such, the sAA convergence factor, which can be computed analytically when the window size $m=1$, can be
%used, as we show in this paper, to quantify the AA-ADMM convergence speed. Our numerical tests will also
%show, as in \cite{desterck2020}, how the spectral properties of the ADMM and optimal sAA-ADMM Jacobians
%can be used to give insight into the mechanism by which nonlinear acceleration methods like sAA and AA
%accelerate ADMM.

The results in \cite{desterck2020} consider sAA with $m = 1$ applied to FPI (\ref{eq:FPI}):
\begin{equation}\label{eq:sAA1}
    \mathbf x_{k+1} = \alpha_0\mathbf q(\mathbf x_k) + \alpha_1 \mathbf q(\mathbf x_{k-1}) = (1+\beta)\mathbf q(\mathbf x_k)-\beta \mathbf q(\mathbf x_{k-1}),
\end{equation}
where $\beta$ remains fixed at all iterations. 
Note that, for $m=1$, this is a stationary version of Nesterov's accelerated gradient descent method
if $\mathbf q(\mathbf x)$ is a gradient descent update.

To study the convergence behaviour and find the optimal choice of $\beta$, we introduce
\begin{equation*}
    \mathbf X_k = 
    \begin{bmatrix}
    \mathbf x_k\\
    \mathbf x_{k-1}
    \end{bmatrix}
\end{equation*}
and write sAA iteration (\ref{eq:sAA1}) as
\begin{equation*}
\mathbf X_{k+1} = 
    \begin{bmatrix}
    \mathbf x_{k+1}\\
    \mathbf x_{k}
    \end{bmatrix}
    =
    \begin{bmatrix}
    (1+\beta)\mathbf q(\mathbf x_k)-\beta \mathbf q(\mathbf x_{k-1})\\
    \mathbf x_k
    \end{bmatrix}
    = \mathbf{\Psi} (\mathbf X_k).
\end{equation*}
% Because $q(\cdot)$ is assumed to be continuously differentiable, we have by Taylor expansion
% \begin{equation*}
%     \mathbf x_{k+1} - \mathbf X^* = \Psi (\mathbf X_k) - \mathbf X^* = \Psi'(\mathbf X^*) (\mathbf X_k - \mathbf X^*) + o(||X_k-X^*||),
% \end{equation*}
Ostrowski's theorem \cite[Theorem 10.1.3]{ortega2000iterative} implies that, when $\mathbf\Psi$ 
has a fixed point $\mathbf X^*$ and is F-differentiable at $\mathbf X^*$, and the spectral radius of $\mathbf\Psi'$ at $\mathbf X^*$ satisfies $\rho(\mathbf\Psi') < 1$, then $\mathbf X^*$ is a point of attraction of the iteration $\mathbf X_{k+1}=\mathbf\Psi(\mathbf X_k)$,
where
\begin{equation*}
    \mathbf\Psi'(\mathbf X^*) = 
    \begin{bmatrix}
    (1+\beta)\mathbf q'(\mathbf x^*) & \ -\beta \mathbf q'(\mathbf x^*)\\
    \mathbf I & \mathbf O
    \end{bmatrix}.
\end{equation*}
% Therefore,
% \begin{equation*}
%     \lim_{k\rightarrow\infty}\frac{||\mathbf x_{k+1}-X^*||}{||X_k-X^*||} \leq ||\Psi'(X^*)|| = ||\Psi'_*||.
% \end{equation*}
% This gives a tight upper bound of the asymptotic convergence factor of sAA(1),
% $\rho_{sAA(1)}=\rho(\Psi'(X^*))$. 
In addition, if $\rho(\mathbf\Psi') > 0$, the iteration will have a root-linear convergence factor that is given by $\rho(\mathbf\Psi')$\cite[Theorem 10.1.4]{ortega2000iterative}. We are interested in finding the optimal asymptotic convergence factor $\rho^*_{sAA(1)}$ of sAA(1)
over all possible choices of $\beta$:
$$
\rho^*_{sAA(1)}=\min_\beta \rho_{sAA(1)}(\beta).
$$
By the properties of the Schur complement, we have that
\begin{align*}
    |\lambda\mathbf I &- \mathbf\Psi'(\mathbf X^*)| = \left|
    \begin{matrix}
    \lambda\mathbf I - (1+\beta)\mathbf q'(\mathbf x^*) & \ \beta \mathbf q'(\mathbf x^*)\\
    -\mathbf I & \lambda \mathbf I
    \end{matrix} 
    \right|\\
    &= 
    \left|
    \lambda\mathbf \, (\lambda\mathbf I - (1+\beta)\mathbf q'(\mathbf x^*)) + \beta \mathbf q'(\mathbf x^*)
    \right|
    = |\lambda^2 \mathbf I -(1+\beta)\lambda \mathbf q'(\mathbf x^*) + \beta\mathbf q'(\mathbf x^*)| = 0
\end{align*}
where $\lambda$ is any eigenvalue of $\mathbf\Psi'(\mathbf X^*)$ and $|\mathbf M|$ means the determinant of matrix $\mathbf M$. Denote the eigenvalues of $\mathbf q'(\mathbf x^*)$ by $\mu$, then we have
\begin{equation}\label{eqn:sAA2_lambda_mu_relation}
    \lambda^2 - (1+\beta)\mu\lambda + \beta\mu = 0.
\end{equation}
Hence, all the eigenvalues of $\mathbf\Psi'(\mathbf X^*)$ are contained in the set 
\begin{equation*}
    \{\lambda : \lambda^2 - (1+\beta)\mu\lambda + \beta\mu = 0, \; \mu\in\sigma(\mathbf q'(\mathbf x^*))\},
\end{equation*}
where $\sigma(\mathbf M)$ means the spectrum of matrix $\mathbf M$. To determine the optimal $\beta$, we only need to find
\begin{equation*}
    \beta^* = \underset{\beta\in\mathbb R}{\arg\min}\max \{|\lambda|: \lambda^2 - (1+\beta)\mu\lambda + \beta\mu = 0, \; \mu\in\sigma(\mathbf q'(\mathbf x^*)) \}.
\end{equation*}
To compute $\beta^*$, we define, for any fixed $\mu$, 
\begin{equation*}
    S_{\mu}(\beta) = \max\{|\lambda|: \lambda^2 - (1+\beta)\mu\lambda + \beta\mu = 0\}.
\end{equation*}
We first assume that the spectrum of $\mathbf q'(\mathbf x^*)$ is real.
Then the following conclusions hold:
\begin{proposition}\label{prop:circle}
Assume $\mu \in \mathbb R$. Any complex eigenvalues $\lambda$ of $\mathbf\Psi'(\mathbf x^*)$ lie on a circle of radius $\left|\frac{\beta}{1+\beta}\right|$ centered at $(\frac{\beta}{1+\beta}, 0)$ in the complex plane.
\begin{proof}
From the relation of $\lambda$ and $\mu$ in (\ref{eqn:sAA2_lambda_mu_relation}), if the roots are complex, i.e. $(1+\beta)^2\mu^2-4\beta\mu < 0$, then
\begin{equation*}
    \lambda, \bar\lambda = \frac{(1+\beta)\mu}{2} \pm i \frac{\sqrt{4\beta\mu-(1+\beta)^2\mu^2}}{2}.
\end{equation*}
Hence, we get
\begin{equation*}
    \lambda\bar\lambda = \beta\mu, \quad \lambda + \bar\lambda = (1+\beta)\mu.
\end{equation*}
Since
\begin{equation*}
    \lambda\bar\lambda - \frac{\beta}{1+\beta}(\lambda+\bar\lambda) + \left(\frac{\beta}{1+\beta}\right)^2 = \left(\frac{\beta}{1+\beta}\right)^2,
\end{equation*}
we have
\begin{equation*}
    \left|\lambda - \frac{\beta}{1+\beta} \right|^2 = \left(\frac{\beta}{1+\beta}\right)^2.
\end{equation*}
This finishes the proof.
\end{proof}
\end{proposition}

\begin{proposition}
\cite[Lemmas 3.1,3.2]{desterck2020} 
When $0 < \mu < 1$, $\min_{\beta}S_{\mu}(\beta) = 1-\sqrt{1-\mu}$, and the optimum is achieved at $\beta_{\mu}^* = \displaystyle \frac{1-\sqrt{1-\mu}}{1+\sqrt{1-\mu}}$.

When $\mu \geq 1$, $\min_{\beta}S_{\mu}(\beta) = \sqrt{\mu}$, and the optimum is achieved at $\beta_{\mu}^*=-1$.

When $\mu < 0$, $\min_{\beta}S_{\mu}(\beta) = \sqrt{1-\mu}-1$, and the optimum is achieved at $\beta_{\mu}^* = 
\displaystyle \frac{1-\sqrt{1-\mu}}{1+\sqrt{1-\mu}}$.
\end{proposition}

From this proposition and the monotonicity of $\min_{\beta} S_{\mu}(\beta)$
over $\mu$ \cite{desterck2020}, 
still for the case the spectrum of $\mathbf q'(\mathbf x^*)$ is real,
we can easily derive the following proposition
where we denote 
\begin{equation*}
    \sigma_{\max}=\max(\sigma(\mathbf q'(\mathbf x^*))), \quad \sigma_{\min}= \min(\sigma(\mathbf q'(\mathbf x^*))).
\end{equation*}

\begin{proposition}\label{prop2}
(Extension of \cite[Theorem 3.4]{desterck2020}.)
When $\sigma(\mathbf q'(\mathbf x^*)) \subset [0, 1)$, the optimal weight is 
$$\beta^* = \frac{1-\sqrt{1-\sigma_{\max}}}{1+\sqrt{1-\sigma_{\max}}},$$
and the optimal convergence factor is $\rho_{sAA(1)}^* = 1- \sqrt{1-\sigma_{\max}}$.

When $\sigma(q'(\mathbf x^*)) \subset (-1, 0]$, the optimal weight is 
$$\beta^* = \frac{1-\sqrt{1-\sigma_{\min}}}{1+\sqrt{1-\sigma_{\min}}},$$ and the optimal convergence factor is $\rho_{sAA(1)}^* = \sqrt{1-\sigma_{\min}}-1$.

When $\sigma({\bf{q}'(\bf{x}^*)})\subset (-1,1)$ and $\sigma_{\rm max}\sigma_{\rm min}<0$, we consider three cases. Define
\begin{equation*}
  \beta_{+} =\frac{1-\sqrt{1-\sigma_{\rm max}}}{1+\sqrt{1-\sigma_{\rm max}}}, \quad   \beta_{-} =\frac{1-\sqrt{1-\sigma_{\rm min}}}{1+\sqrt{1-\sigma_{\rm min}}}.
\end{equation*}
\begin{itemize}
\item[(a)] If $\sigma_{\rm max} =|\sigma_{\rm min}|$, then the  optimal weight is $\beta^*=0$ and $\rho^*_{sAA(1)}=\sigma_{\rm max}$.

\item[(b)] If  $\sigma_{\rm max} >|\sigma_{\rm min}|$, there are two subcases:\\
(b1) If $\displaystyle \frac{-(1+\beta_+)\sigma_{\rm min}+\sqrt{(1+\beta_+)^2\sigma_{\rm min}^2-4\beta_+\sigma_{\rm min}}}{2}\leq 1-\sqrt{1-\sigma_{\rm max}}$, then
\begin{equation*}
  \beta^*=\beta_+,\quad \rho^*_{sAA(1)} = 1-\sqrt{1-\sigma_{\rm max}}.
\end{equation*}

(b2) If $\displaystyle \frac{-(1+\beta_+)\sigma_{\rm min}+\sqrt{(1+\beta_+)^2\sigma_{\rm min}^2-4\beta_+\sigma_{\rm min}}}{2}> 1-\sqrt{1-\sigma_{\rm max}}$, then
the optimal $\beta^*$ is obtained by solving
\begin{equation*}
  \frac{-(1+\beta)\sigma_{\rm min}+\sqrt{(1+\beta)^2\sigma_{\rm min}^2-4\beta\sigma_{\rm min}}}{2}=\frac{(1+\beta)\sigma_{\rm max}+\sqrt{(1+\beta)^2\sigma_{\rm max}^2-4\beta\sigma_{\rm max}}}{2},
\end{equation*}
which gives
\begin{equation*}
\beta^* = \frac{(m_+-\sqrt{m_+^2-4})^2}{4}, \quad \text{where}\quad m_+ =\frac{\sigma_{\rm \max}-\sigma_{\rm min}}{\sqrt{-2\sigma_{\rm max}\sigma_{\rm min}(\sigma_{\rm max}+\sigma_{\rm min})}},
\end{equation*}
and the corresponding optimal convergence factor is
\begin{equation*}
  \rho^*_{sAA(1)} = \frac{(1+\beta^*)\sigma_{\rm max}+\sqrt{(1+\beta^*)^2\sigma_{\rm max}^2-4\beta^*\sigma_{\rm max}}}{2} > 1-\sqrt{1-\sigma_{\rm max}}.
\end{equation*}

\item[(c)] If  $\sigma_{\rm max} <|\sigma_{\rm min}|$, there are two subcases:

(c1) If $\displaystyle \frac{(1+\beta_-)\sigma_{\rm max}+\sqrt{(1+\beta_-)^2\sigma_{\rm max}^2-4\beta_-\sigma_{\rm max}}}{2}\leq \sqrt{1-\sigma_{\rm min}}-1$, then
\begin{equation*}
  \beta^*=\beta_-,\quad \rho^*_{sAA(1)} = \sqrt{1-\sigma_{\rm min}}-1.
\end{equation*}

(c2) If $\displaystyle \frac{(1+\beta_-)\sigma_{\rm max}+\sqrt{(1+\beta_-)^2\sigma_{\rm max}^2-4\beta_-\sigma_{\rm max}}}{2}> \sqrt{1-\sigma_{\rm min}}-1$, then
the optimal $\beta^*$ is obtained by solving
\begin{equation*}
  \frac{-(1+\beta)\sigma_{\rm min}+\sqrt{(1+\beta)^2\sigma_{\rm min}^2-4\beta\sigma_{\rm min}}}{2}=\frac{(1+\beta)\sigma_{\rm max}+\sqrt{(1+\beta)^2\sigma_{\rm max}^2-4\beta\sigma_{\rm max}}}{2},
\end{equation*}
which gives
\begin{equation*}
\beta^* = -\frac{(\sqrt{m_{-}^2+4}-m_{-})^2}{4}, \quad \text{where}\quad m_{-} =\frac{\sigma_{\rm \max}-\sigma_{\rm min}}{\sqrt{2\sigma_{\rm max}\sigma_{\rm min}(\sigma_{\rm max}+\sigma_{\rm min})}},
\end{equation*}
and the corresponding optimal convergence factor is
\begin{equation*}
  \rho^*_{sAA(1)} = \frac{(1+\beta^*)\sigma_{\rm max}+\sqrt{(1+\beta^*)^2\sigma_{\rm max}^2-4\beta^*\sigma_{\rm max}}}{2}>\sqrt{1-\sigma_{\rm min}}-1.
\end{equation*}

\end{itemize}

%When $\sigma(\mathbf q'(\mathbf x^*)) \subset (-1, 1)$ and $\sigma_{\max}\sigma_{\min} < 0$, the optimal weight is
%\begin{equation*}
%    \beta^* = 
%    \begin{cases}
%        \displaystyle \frac{1-\sqrt{1-\sigma_{\max}}}{1+\sqrt{1-\sigma_{\max}}} & \textrm{if } \ \sigma_{\max} > |\sigma_{\min}|, \\
%        \\
%        \displaystyle \frac{1-\sqrt{1-\sigma_{\min}}}{1+\sqrt{1-\sigma_{\min}}} & \textrm{if } \ \sigma_{\max} \le |\sigma_{\min}|,
%    \end{cases}
%\end{equation*}
%and the optimal convergence factor is
%\begin{equation*}
%    \rho^*_{sAA(1)} = 
%    \begin{cases}
%        1-\sqrt{1-\sigma_{\max}} & \textrm{if } \ \sigma_{\max} > |\sigma_{\min}|, \\
%        \sqrt{1-\sigma_{\min}}-1 & \textrm{if } \ \sigma_{\max} \le |\sigma_{\min}|.
%    \end{cases}
%\end{equation*}
\end{proposition}

\begin{remark}
The result for $\sigma_{\rm max}\sigma_{\rm \min}<0$ is an extension of Theorem 3.4 in \cite{desterck2020}, and follows directly from the proof there. This case does not occur in the test problems we consider in this paper, but we include it for completeness since it may arise in other applications.
\end{remark}

If the spectrum of $\mathbf q'(\mathbf x^*)$ is complex, the following result can be used:
\begin{proposition}\label{prop3}
\cite{desterck2020} Let the spectral radius of $\mathbf q'(\mathbf x^*)$ be $\rho_{q'}^*$ and assume $\rho_{q'}^*<1$.
If there exists a real eigenvalue $\mu$ of $\mathbf q'(\mathbf x^*)$ such that $\rho_{ q'}^* = \mu$, then the optimal asymptotic convergence rate of sAA(1), $\rho^*_{sAA(1)}$, is bounded below by
\begin{equation*}
    \rho^*_{sAA(1)} \geq 1-\sqrt{1-\rho^*_{q'}},
\end{equation*}
and if the equality holds,
\begin{equation*}
    \beta^* = \frac{1-\sqrt{1-\rho^*_{q'}}}{1+\sqrt{1-\rho^*_{q'}}}.
\end{equation*}
\end{proposition}

Propositions \ref{prop2} and \ref{prop3} allow us to compute the optimal sAA(1) coefficient $\beta^*$ and the optimal
asymptotic convergence factor, $\rho^*_{sAA}$, (or a lower bound) when $\mathbf q'(\mathbf x^*)$ is known.
%In \cite{desterck2020}, the case of $q'(\mathbf x^*)$ with complex spectrum is also discussed, but we do not need this
%case for the numerical results discussed in this paper.
Also, \cite{desterck2020} explains how optimal sAA weights and convergence factors $\rho^*_{sAA}$ can be determined for sAA with $m\ge 2$ by optimization, since analytical results are not known in this case. For example, for the case when $m = 2$, the sAA(2) iteration is
\begin{equation}
    \mathbf x_{k+1} = (1+\beta_1+\beta_2)\mathbf q(x_k) - \beta_1\mathbf q(x_{k-1}) - \beta_2\mathbf q(x_{k-2}).
     \label{eq:sAA(2)}
\end{equation}
We compute the optimal $\beta_1^*$ and $\beta_2^*$ from
\begin{equation*}
    \{\beta_1^*,\beta_2^*\} = \underset{\beta_1,\beta_2\in\mathbb R}{\arg\min}\max_{\lambda} \{|\lambda|: \lambda^3 - (1+\beta_1+\beta_2)\mu\lambda^2 + \beta_1\mu\lambda + \beta_2\mu = 0, \; \mu\in\sigma(\mathbf q'(\mathbf x^*)) \},
\end{equation*}
which can be solved, for example, by brute-force search.
%\qHa{Dawei, I will ask you later to add the optimal asymptotic convergence factor, $\rho^*_{sAA}$, to \autoref{prop2}.}
% \begin{figure}[H]
% \centering
% \includegraphics[scale=0.44]{speedup.jpg}
% \caption{Speedup of AA applied to FPI for different $\mu$}
% \label{fig:aa_speedup}
% \end{figure}

%================================
\section{Acceleration of ADMM by optimal stationary AA and comparison with non-stationary AA}
\label{sec:results}
%================================
In this section we present results analyzing how the optimal convergence factor
of the stationary AA method with window size $m=1$, as computed from \autoref{prop2} and \autoref{prop3},
improves the ADMM convergence speed.
We also consider acceleration by stationary AA with window sizes $m=2$
and $m=3$, where the optimal sAA coefficients are determined by optimization.
We consider a variety of ADMM examples that include linear and nonlinear
cases, smooth and non-smooth cases, and cases with real and complex Jacobian spectrum. 
We investigate the spectra of the ADMM and optimal sAA-ADMM Jacobians
to explain the convergence acceleration
and compare numerically with the asymptotic convergence speed
of ADMM accelerated by non-stationary AA with finite window size.

In all numerical experiments, we use a zero initial guess
unless stated otherwise, and no parameter tuning is applied. 
For the sAA(2) iteration (\ref{eq:sAA(2)}), we approximate the optimal $\beta_1^*$ and $\beta_2^*$
using brute-force search in the range of $[-1, 1]$ with step size $0.05$.
Similarly, we also include some simulation results for sAA(3), using
a brute-force search technique to approximate the optimal $\beta$'s.
Since this approach is expensive for $m=3$, we only report results
for a selection of our test problems.

%As mentioned before, the analysis of sAA with $m \geq 2$ requires solving for the roots of high order polynomials, which is hard. As a result, we are not able to give analytical results of the optimal $\beta$'s for the cases when $m \geq 2$. Instead, we use brute-force search technique to find the optimal $\beta$'s within certain ranges for each specific problem. This naive technique becomes very computationally expensive when $m$ becomes larger. Therefore, we only give numerical results of using sAA(3) for ridge regression and logistic regression problems, each representing linear and nonlinear problem respectively, to emphasize our conjecture that the asymptotic convergence factor of sAA($m$) gives a useful prediction of that of AA($m$).

%-------------------------------------------------------
\subsection{Ridge regression (see, e.g., \cite{boyd2011distributed}; linear and smooth problem)}
%-------------------------------------------------------
%,,,,,,,,,,,,,,,,,,,,,,,,,,,,,,,,,,,,,,,,,,,,,,,,,,,,,,,,,,,,,,,,,,,,,,,,,,,,,,
\subsubsection{Problem description}
%,,,,,,,,,,,,,,,,,,,,,,,,,,,,,,,,,,,,,,,,,,,,,,,,,,,,,,,,,,,,,,,,,,,,,,,,,,,,,,
The $l_2$-regularized least squares problem, also called ridge regression, is a common technique in machine learning that reduces model complexity and prevents over-fitting. The optimization problem is 
\begin{equation*}
    \min_{\mathbf x}\frac{1}{2}||\mathbf A \mathbf x - \mathbf b||_2^2 + \lambda ||\mathbf x||_2^2,
\end{equation*}
where $(\mathbf A, \mathbf b)\in\mathbb R^{m\times n}\times\mathbb R^m$ is the training set, and $\lambda > 0$ is a regularization parameter.

To use the ADMM method, we write this problem as
\begin{align}\label{eq:ridge}
    &\min_{\mathbf x, \mathbf z}  \frac{1}{2}||\mathbf A \mathbf x - b||_2^2 + \lambda ||\mathbf z||_2^2, \\
    &\text{s.t. } \mathbf x - \mathbf z  = 0. \nonumber
\end{align}
The scaled augmented Lagrangian is
\begin{align*}
    L_{\rho}(\mathbf x, \mathbf z, \mathbf u) = \frac{1}{2}||\mathbf A \mathbf x -\mathbf b||_2^2 + \lambda ||\mathbf z||_2^2 + \frac{\rho}{2} ||\mathbf x - \mathbf z + \mathbf u||_2^2 - \frac{\rho}{2}||\mathbf u||_2^2.
\end{align*}
The ADMM steps for this problem are:
\begin{equation*}
\begin{cases}
    \mathbf x_{k+1} = \text{argmin}_{\mathbf x} \frac{1}{2}||\mathbf A \mathbf x - \mathbf b||_2^2 + \frac{\rho}{2} ||\mathbf x - \mathbf z_k + \mathbf u_k||_2^2,\\
    \mathbf z_{k+1} = \text{argmin}_{\mathbf z}\lambda ||\mathbf z||_2^2 + \frac{\rho}{2} ||\mathbf x_{k+1}  + \mathbf u_k - \mathbf z||_2^2, \\
    \mathbf u_{k+1} = \mathbf u_k + \mathbf x_{k+1} - \mathbf z_{k+1},
\end{cases}
\end{equation*}
which gives
\begin{equation*}
\begin{cases}
    \mathbf x_{k+1} = (\mathbf A^T\mathbf A+\rho \mathbf I)^{-1}\left( \mathbf A^T\mathbf b + \rho(\mathbf z_k - \mathbf u_k) \right)\\
    \mathbf z_{k+1} = \frac{\rho}{2\lambda+\rho}(\mathbf x_{k+1}+\mathbf u_k) \\
    \mathbf u_{k+1} = \mathbf u_k + \mathbf x_{k+1} - \mathbf z_{k+1}.
\end{cases}
\end{equation*}
%where $x^{k+1}$ is the proximal operator of the $l_2$ norm which can be evaluated from a least squares problem
%\begin{equation}\label{eqn:prox_l2_x_update}
%    \mathbf x_{k+1} = \text{argmin}_{\mathbf x}\left|\left| 
%    \begin{bmatrix}
%    A \\
%    \sqrt{\rho}I
%    \end{bmatrix} x -
%    \begin{bmatrix}
%    b \\
%    \sqrt{\rho} (z_k - u_k)
%    \end{bmatrix}
%    \right|\right|_2^2.
%\end{equation}
Since $\mathbf u_{k+1}$ can be explicitly obtained from $\mathbf z_{k+1}$,
\begin{equation*}
    \mathbf u_{k+1} = \frac{2\lambda}{\rho} \mathbf z_{k+1},
\end{equation*}
we can write one iteration of ADMM as a fixed-point update of variable $\mathbf z$, $\bar{\mathbf z}_{k+1} = \mathbf q(\mathbf z_k)$, where
\begin{align*}
    \mathbf q(\mathbf z_k) &= \left[ \frac{\rho(\rho-2\lambda)}{\rho+2\lambda}(\mathbf A^T\mathbf A+\rho \mathbf I)^{-1}+\frac{2\lambda}{\rho+2\lambda}\mathbf I \right]\mathbf z_k + \frac{\rho}{\rho+2\lambda}(\mathbf A^T\mathbf A+\rho \mathbf I)^{-1}\mathbf A^T\mathbf b \\
    &= \mathbf M\mathbf z_k+\hat{\mathbf b}.
\end{align*}
Problem (\ref{eq:ridge}) has the closed-form exact solution
\begin{equation*}
    \mathbf x^* = \mathbf z^* = (\mathbf A^T\mathbf A+2\lambda \mathbf I)^{-1}\mathbf A^T \mathbf b,
\end{equation*}
and solving it by ADMM is not of practical interest.
Still, convergence acceleration of ADMM for this problem is interesting for our purposes, since
it illustrates our approach and results in the most simple linear and smooth setting,
and will be followed by increasingly complex nonlinear and non-smooth problems in
our further examples. 
The ADMM update is simply a stationary linear iteration. Therefore, $\mathbf q' = \mathbf M$ is independent of $\mathbf z$. To determine the optimal sAA acceleration, we can analyze the spectrum of matrix $\mathbf M$ to pick the optimal $\beta^*$. 

%,,,,,,,,,,,,,,,,,,,,,,,,,,,,,,,,,,,,,,,,,,,,,,,,,,,,,,,,,,,,,,,,,,,,,,,,,,,,,,
\subsubsection{Parameters for test problem}
%,,,,,,,,,,,,,,,,,,,,,,,,,,,,,,,,,,,,,,,,,,,,,,,,,,,,,,,,,,,,,,,,,,,,,,,,,,,,,,
We implement our algorithms on a randomly generated sparse matrix of size $m\times n = 150\times300$ with density 0.001
%(ratio of the number of nonzeros) 
sampled from the standard normal distribution. The $\mathbf b$ vector is sampled from the standard normal distribution. The regularization parameter is chosen as $\lambda = 1$, and we pick the penalty parameter $\rho = 10$.

\begin{figure}[t]
\centering
\includegraphics[scale=0.65]{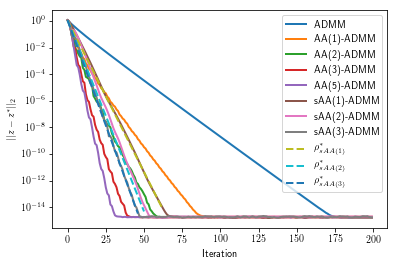}
\includegraphics[scale=0.65]{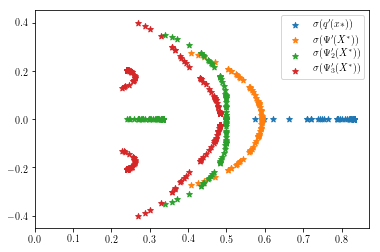}
\caption{Ridge regression. (top) Comparison of error reduction using ADMM, AA($m$)-ADMM and sAA($m$)-ADMM.
(bottom) Spectrum of $\mathbf q'$ of ADMM, $\mathbf\Psi'$ of sAA(1)-ADMM, and $\mathbf\Psi_2'$ and $\mathbf\Psi_3'$ of sAA(2)-ADMM and sAA(3)-ADMM.}
\label{fig:ridge_e}
\end{figure}

%\begin{figure}[H]
%\centering
%\caption{Ridge regression. (top) Comparison of error reduction using ADMM, sAA($m$)-ADMM and AA($m$)-ADMM.
%(bottom) Spectrum of ADMM iteration matrix $\mathbf q'$ and sAA($m$)-ADMM iteration matrix $\mathbf\Psi'_m$.}
%\label{fig:ridge_e_sAA3}
%\end{figure}

%,,,,,,,,,,,,,,,,,,,,,,,,,,,,,,,,,,,,,,,,,,,,,,,,,,,,,,,,,,,,,,,,,,,,,,,,,,,,,,
\subsubsection{Convergence results} 
%,,,,,,,,,,,,,,,,,,,,,,,,,,,,,,,,,,,,,,,,,,,,,,,,,,,,,,,,,,,,,,,,,,,,,,,,,,,,,,

%\begin{figure}[H]
%\centering
%\centering
%\includegraphics[scale=0.65]{figures/eigs_ridge.png}
%\caption{Ridge regression: }
%\label{fig:eigs_ridge}
%\end{figure}

We obtain convergence plots for the error $\|\mathbf z-\mathbf z^*\|_2$ as shown in \autoref{fig:ridge_e} (top).
We see that ADMM converges linearly.
The convergence factor of ADMM is substantially improved by the AA-based methods.
AA(2) and AA(3) converge slightly faster than AA(1), and sAA(1) converges with similar
asymptotic speed.

The convergence improvement of the AA-ADMM methods over ADMM can be understood in terms
of spectral properties as follows.
\autoref{fig:ridge_e} (bottom) shows the spectrum of the ADMM iteration matrix $\mathbf M$, $\sigma(\mathbf M)\in(0,1)$.
The spectrum is real since $\mathbf M$ is symmetric, and the spectral radius $\rho^*_{q'} = 0.833.$
Therefore, according to \autoref{prop2}, the optimal $\beta$ for sAA(1) is
\begin{equation*}
    \beta^* = \frac{1-\sqrt{1-\rho^*_{q'}}}{1+\sqrt{1-\rho^*_{q'}}} = 0.420.
\end{equation*}
The corresponding optimal sAA(1) linear convergence factor is 
\begin{equation*}
    \rho^*_{sAA(1)-ADMM} = \rho(\mathbf \Psi') = 1 - \sqrt{1-\rho^*_{q'}} = 0.592<0.833.
\end{equation*}
The approximately optimal $\beta_1^*$ and $\beta_2^*$ for sAA(2) are
\begin{equation*}
    \beta_1^* = 0.70, \; \beta_2^* = -0.10,
\end{equation*}
with sAA(2) linear convergence factor
\begin{equation*}
    \rho^*_{sAA(2)-ADMM} = \rho(\mathbf\Psi_2') = 0.516.
\end{equation*}
Similarly, for sAA(3), we obtain
\begin{equation*}
    \beta_1^* = 0.955, \; \beta_2^* = -0.250, \; \beta_3^* = 0.028,
\end{equation*}
with sAA(3) linear convergence factor
\begin{equation*}
    \rho^*_{sAA(3)-ADMM} = \rho(\mathbf\Psi_3') = 0.4837<\rho(\mathbf\Psi_2')<\rho(\mathbf\Psi')<\rho^*_{q'}.
\end{equation*}

\autoref{fig:ridge_e} (bottom) also shows the spectrum of the sAA(1)-ADMM iteration matrix, $\mathbf\Psi'$, and of the sAA(2)-ADMM and sAA(3)-ADMM iteration matrices, $\mathbf\Psi_2'$ and $\mathbf\Psi_3'$.
The acceleration methods spread the ADMM spectrum out in the complex plane
in a way that strongly reduces the asymptotic convergence factor: 
e.g., $\rho(\mathbf\Psi')$ is much smaller than $\rho^*_{q'}$.
Note that stationary iterative method (\ref{eq:sAA1}) maps part of the nonnegative real spectrum of $\mathbf q'(\mathbf x^*)$ to a circle, according to \autoref{prop:circle}.
As seen in \autoref{fig:ridge_e} (top), the optimal sAA(1)-ADMM factor, $\rho^*_{sAA(1)-ADMM}$, provides a
useful prediction of the convergence factors of the AA-ADMM methods.
The convergence speed of sAA(1)-ADMM matches the theoretical prediction of $\rho^*_{sAA(1)-ADMM}$.

%Finally, \autoref{fig:ridge_r} shows the reduction in the combined residual for the different methods.
%\begin{figure}
%\centering
%\includegraphics[scale=0.65]{figures/ridge_r.png}
%\caption{Ridge regression: comparison of residual reduction using ADMM, AA(1)-ADMM and sAA(1)-ADMM}
%\label{fig:ridge_r}
%\end{figure}

%-------------------------------------------------------
\subsection{Regularized logistic regression (see, e.g., \cite{boyd2011distributed}; nonlinear and smooth problem)}\label{sec:reglog}
%-------------------------------------------------------
%,,,,,,,,,,,,,,,,,,,,,,,,,,,,,,,,,,,,,,,,,,,,,,,,,,,,,,,,,,,,,,,,,,,,,,,,,,,,,,
\subsubsection{Problem description}
%,,,,,,,,,,,,,,,,,,,,,,,,,,,,,,,,,,,,,,,,,,,,,,,,,,,,,,,,,,,,,,,,,,,,,,,,,,,,,,
We consider a simple logistic regression model in this section. The objective function of the regularized logistic regression model is
\begin{equation*}
    \min_{\mathbf x} \frac{1}{m}\sum_{i=1}^m\log(1 + \exp(-y_i(\mathbf a_i^T\mathbf w + \mathbf c)) + \lambda ||\mathbf x||_2^2,
\end{equation*}
where 
\begin{equation*}
    A = 
    \begin{bmatrix}
    \mathbf a_1^T\\
    \vdots\\
    \mathbf a_m^T
    \end{bmatrix} \in \mathbb R^{m\times n},
\end{equation*}
are $m$ data samples, $y_1, \cdots, y_m$ are the corresponding labels, and
\begin{equation*}
    \mathbf x = 
    \begin{bmatrix}
    \mathbf c\\
    \mathbf w
    \end{bmatrix},
    \quad \mathbf w\in \mathbb R^{n}, \quad \mathbf c \in \mathbb R,
\end{equation*}
are the linear combination coefficients and bias to be optimized. To apply ADMM, we write this problem as
\begin{align*}
    \min_{\mathbf x, \mathbf z}\; &\frac{1}{m}\sum_{i=1}^m\log(1 + \exp(-y_i(\mathbf a_i^T\mathbf w + \mathbf c))) + \lambda ||\mathbf z||_2^2,\\
    \text{s.t.}\;& \mathbf x - \mathbf z = 0.
\end{align*}
This gives the augmented Lagrangian
\begin{equation*}
    L(\mathbf x, \mathbf z, \mathbf u,\rho) = \frac{1}{m}\sum_{i=1}^m\log(1 + \exp(-y_i(\mathbf a_i^T\mathbf w + \mathbf c))) + \lambda ||\mathbf z||_2^2 + \frac{\rho}{2}||\mathbf x - \mathbf z+\mathbf u||_2^2 - \frac{\rho}{2}||\mathbf u||_2^2.
\end{equation*}
Hence, we get the ADMM steps
\begin{equation*}
\begin{cases}
    \mathbf x_{k+1} = \text{argmin}_{\mathbf x} \frac{1}{m}\sum_{i=1}^m\log(1 + \exp(-y_i(\mathbf a_i^T\mathbf w + \mathbf c))) + \frac{\rho}{2}||\mathbf x - \mathbf z_k+\mathbf u_k||_2^2,\\
    \mathbf z_{k+1} = \text{argmin}_{\mathbf z} \lambda ||\mathbf z||_2^2 + \frac{\rho}{2}||\mathbf x_{k+1}-\mathbf z+\mathbf u_k||_2^2, \\
    \mathbf u_{k+1} = \mathbf u_k + \mathbf x_{k+1} - \mathbf z_{k+1}.
\end{cases}
\end{equation*}
To solve for $\mathbf x_{k+1}$, we use Newton's method.
%can solve an $l_2$ regularized logistic regression problem. Or we can solve the gradient equation using for example 
%\begin{equation*}
%    \frac{1}{m}\sum_{i=1}^{m}\frac{e^{-y_i(a_i^Tw+c)}}{1+e^{-y_i(a_i^Tw+c)}}(-y_ia_i^T) + \rho ( \mathbf x - \mathbf z_k+\mathbf u_k) = \frac{1}{m}C^T\frac{[e^{Cx}]}{[1+e^{Cx}]}+ \rho ( \mathbf x - \mathbf z_k+\mathbf u_k) = 0,
%\end{equation*}
%where $C = [-y,\; -\text{diag}(\mathbf y)A]$. We can derive the Hessian 
%\begin{equation*}
%    H = \frac{1}{m}C^T\text{diag}\left(\frac{[e^{Cx}]}{[1+e^{Cx}]^2}\right)C+ \rho I
%\end{equation*}
%To compute the spectral radius of the ADMM iteration function $q'(\mathbf z)$, we use finite dirrerence approach to approximate the Jacobian matrix at the solution.

\begin{figure}[t]
\centering
\includegraphics[scale=0.65]{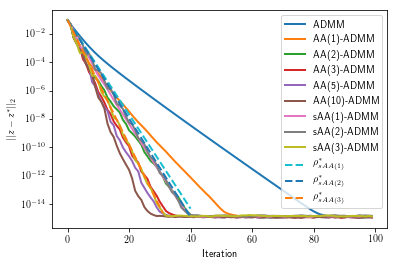}
\includegraphics[scale=0.65]{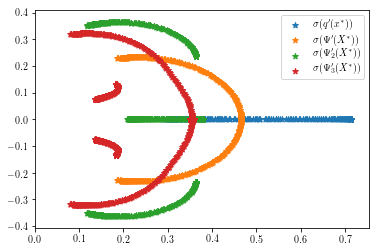}
\caption{$l_2$-regularized logistic regression.
(top) Comparison of error reduction using ADMM, AA($m$)-ADMM and sAA($m$)-ADMM.
(bottom) Spectrum of $\mathbf q'$ of ADMM, $\mathbf\Psi'$ of sAA(1)-ADMM, and $\mathbf\Psi_2'$ and $\mathbf\Psi_3'$ of sAA(2)-ADMM and sAA(3)-ADMM.}
\label{fig:reglog_e}
\end{figure}

%\begin{figure}
%\centering
%\includegraphics[scale=0.65]{figures/logreg_e_sAA3.png}
%\includegraphics[scale=0.65]{figures/eigs_logreg_sAA3.png}
%\caption{Regularized logistic regression. (top) Comparison of error reduction using ADMM, sAA($m$)-ADMM and AA($m$)-ADMM.
%(bottom) Spectrum of ADMM iteration matrix $\mathbf q'$ and sAA($m$)-ADMM iteration matrix $\mathbf\Psi'_m$.}
%\label{fig:logreg_e_sAA3}
%\end{figure}

%\begin{figure}
%\centering
%\caption{$l_2$-regularized logistic regression: .}
%\label{fig:eigs_reglog}
%\end{figure}

%,,,,,,,,,,,,,,,,,,,,,,,,,,,,,,,,,,,,,,,,,,,,,,,,,,,,,,,,,,,,,,,,,,,,,,,,,,,,,,
\subsubsection{Parameters for the test problem}
%,,,,,,,,,,,,,,,,,,,,,,,,,,,,,,,,,,,,,,,,,,,,,,,,,,,,,,,,,,,,,,,,,,,,,,,,,,,,,,
For this problem, we applied our algorithms to the Madelon data set from the UCI machine learning repository\footnote{https://archive.ics.uci.edu/ml/datasets/Madelon}. To reduce the amount of computation, we only used a portion of the features and examples. The regularization parameter is $\lambda = 2$, and the augmented Lagrangian penalty parameter is $\rho = 10$.

%,,,,,,,,,,,,,,,,,,,,,,,,,,,,,,,,,,,,,,,,,,,,,,,,,,,,,,,,,,,,,,,,,,,,,,,,,,,,,,
\subsubsection{Convergence results}
%,,,,,,,,,,,,,,,,,,,,,,,,,,,,,,,,,,,,,,,,,,,,,,,,,,,,,,,,,,,,,,,,,,,,,,,,,,,,,,
Since the FPI representation of ADMM for solving the regularized logistic regression problem is nonlinear, we are now not able to find an explicit expression for $\mathbf z_{k+1}= \mathbf q(\mathbf z_k)$ like before. To determine the spectrum of $\mathbf q'(\mathbf z^*)$, we use the first-order finite difference method with step size $h = 1\times10^{-4}$ to approximate $\mathbf q'(\mathbf z^*)$ at the approximate true solution solved to $10^{-16}$ accuracy.   

\autoref{fig:reglog_e} (top) compares the error norm reduction when using ADMM, AA($m$)-ADMM and sAA($m$)-ADMM.
The convergence acceleration seen in the figure can be explained based on the spectra in
\autoref{fig:reglog_e} (bottom).
The spectrum of $\mathbf q'(\mathbf z^*)$ has asymptotic convergence factor $\rho^*_{q'} = 0.714$.
We can choose the optimal $\beta^*$ the same way as in the ridge regression problem:
\begin{equation*}
    \beta^* = \frac{1-\sqrt{1-\rho^*_{q'}}}{1+\sqrt{1-\rho^*_{q'}}} = 0.303.
\end{equation*}
The corresponding optimal sAA(1)-ADMM linear convergence factor is 
\begin{equation*}
    \rho^*_{sAA(1)-ADMM} = \rho(\mathbf\Psi') = 1 - \sqrt{1-\rho^*_{q'}} = 0.465 < (\rho^*_{q'})^2.
\end{equation*}
The approximately optimal $\beta_1^*$ and $\beta_2^*$ for sAA(2) are
\begin{equation*}
    \beta_1^* = 0.65, \; \beta_2^* = -0.10,
\end{equation*}
with sAA(2) linear convergence factor
\begin{equation*}
    \rho^*_{sAA(2)-ADMM} = \rho(\mathbf\Psi_2') = 0.450.
\end{equation*}
Similarly, for sAA(3), we obtain
\begin{equation*}
    \beta_1^* = 0.61, \; \beta_2^* = -0.115, \; \beta_3^* = 0.009,
\end{equation*}
with sAA(3) linear convergence factor
\begin{equation*}
    \rho^*_{sAA(3)-ADMM} = \rho(\mathbf\Psi_3') = 0.364.
\end{equation*}
\autoref{fig:reglog_e} (top) shows that $\rho^*_{sAA(m)-ADMM}$ is a useful prediction for the
convergence factors of the AA-accelerated ADMM methods.

\subsection{Total variation (see, e.g., \cite{boyd2011distributed}; nonlinear and nonsmooth problem, complex spectrum)}\label{sec:totvar}
%-------------------------------------------------------
%,,,,,,,,,,,,,,,,,,,,,,,,,,,,,,,,,,,,,,,,,,,,,,,,,,,,,,,,,,,,,,,,,,,,,,,,,,,,,,
\subsubsection{Problem description}
%,,,,,,,,,,,,,,,,,,,,,,,,,,,,,,,,,,,,,,,,,,,,,,,,,,,,,,,,,,,,,,,,,,,,,,,,,,,,,,
The total variation model is a widely used method
for applications like image denoising. The optimization problem is
\begin{equation*}
    \min_{\mathbf x} \frac{1}{2}||\mathbf y - \mathbf x||_2^2 + \alpha||\mathbf D\mathbf x||_1^2,
\end{equation*}
where $\mathbf x \in \mathbb R^n$ is the variable, $\mathbf y \in \mathbb R^n$ is the problem data (e.g. image pixel values), $\alpha > 0$ is a smoothing parameter, and $\mathbf D \in \mathbb R^{(n-1)\times n}$ is the difference operator 
\begin{equation*}
\mathbf D =
    \begin{bmatrix}
    -1 & 1 \\
      & -1  & 1 \\
      &    & \ddots & \ddots \\
      &    &        & -1      & 1
    \end{bmatrix}.
\end{equation*}
To use ADMM, we write this problem as
\begin{align*}
    \min_{\mathbf x, \mathbf z}\; &\frac{1}{2}||\mathbf y - \mathbf x||_2^2 + \alpha||\mathbf z||_1^2,\\
    \text{s.t.}\; &\mathbf D\mathbf x-\mathbf z = 0.
\end{align*}
The augmented Lagrangian is
\begin{equation*}
    L_{\rho}(\mathbf x, \mathbf z, \mathbf u) = \frac{1}{2}||\mathbf y - \mathbf x||_2^2 + \alpha||\mathbf z||_1^2 + \frac{\rho}{2}||\mathbf D\mathbf x - \mathbf z + \mathbf u||_2^2 - \frac{\rho}{2}||\mathbf u||_2^2.
\end{equation*}
The ADMM steps for this problem are:
\begin{equation*}
\begin{cases}
    \mathbf x_{k+1} = \text{argmin}_{\mathbf x} \frac{1}{2}||\mathbf y - \mathbf x||_2^2 + \frac{\rho}{2}||\mathbf D\mathbf x - \mathbf z_k + \mathbf u_k||_2^2,\\
    \mathbf z_{k+1} = \text{argmin}_{\mathbf z} \alpha||\mathbf z||_1^2 + \frac{\rho}{2} ||\mathbf D\mathbf x_{k+1}  +\mathbf u_k - \mathbf z||_2^2, \\
    \mathbf u_{k+1} = \mathbf u_k + \mathbf D\mathbf x_{k+1} - \mathbf z_{k+1},
\end{cases}
\end{equation*}
where $\mathbf x_{k+1}$ is the proximal operator of the $l_2$ norm which can be evaluated from a least squares problem as before, 
\begin{equation*}
    \mathbf x_{k+1} = \text{argmin}_{\mathbf x}\left|\left| 
    \begin{bmatrix}
    \mathbf D \\
    \frac{1}{\sqrt{\rho}}\mathbf I
    \end{bmatrix} x -
    \begin{bmatrix}
    \mathbf z_k - \mathbf u_k \\
    \frac{1}{\sqrt{\rho}} \mathbf y
    \end{bmatrix}
    \right|\right|_2^2,
\end{equation*}
and $\mathbf z^{k+1}$ is just the proximal operator of the $l_1$-norm,
\begin{equation*}
    \mathbf z_{k+1} = \text{prox}_{\frac{\alpha}{\rho}||\cdot||_1}(\mathbf D\mathbf x_{k+1}+\mathbf u_k).
\end{equation*}

\begin{figure}
\centering
\includegraphics[scale=0.65]{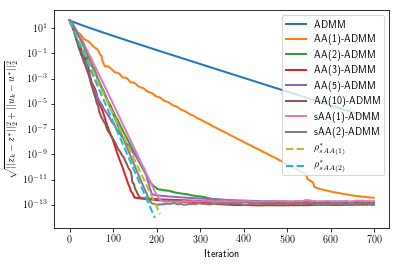}
\includegraphics[scale=0.65]{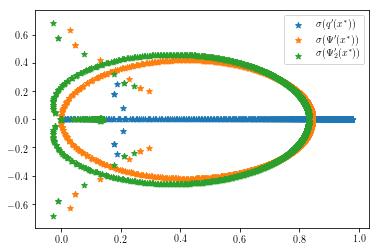}
\caption{Total variation. (top) Comparison of error reduction using ADMM, AA($m$)-ADMM and sAA($m$)-ADMM.
(bottom) Spectrum of ADMM iteration matrix $\mathbf q'$ and sAA($1$)-ADMM iteration matrix $\mathbf\Psi'$.}
\label{fig:total_variation_e}
\end{figure}

%\begin{figure}
%\centering
%\centering
%\caption{Total variation: }
%\label{fig:eigs_total_variation}
%\end{figure}

%,,,,,,,,,,,,,,,,,,,,,,,,,,,,,,,,,,,,,,,,,,,,,,,,,,,,,,,,,,,,,,,,,,,,,,,,,,,,,,
\subsubsection{Parameters for the test problem}
%,,,,,,,,,,,,,,,,,,,,,,,,,,,,,,,,,,,,,,,,,,,,,,,,,,,,,,,,,,,,,,,,,,,,,,,,,,,,,,
We test our algorithms on randomly generated data $\mathbf y$ of size $1000$ sampled from the standard normal distribution. The smoothing parameter is $\alpha = 0.001\cdot||\mathbf y||_{\infty}$. For the penalty parameter, we use $\rho = 10$.

%,,,,,,,,,,,,,,,,,,,,,,,,,,,,,,,,,,,,,,,,,,,,,,,,,,,,,,,,,,,,,,,,,,,,,,,,,,,,,,
\subsubsection{Convergence results}
%,,,,,,,,,,,,,,,,,,,,,,,,,,,,,,,,,,,,,,,,,,,,,,,,,,,,,,,,,,,,,,,,,,,,,,,,,,,,,,
We use the first-order finite difference method with step size $h = 1\times 10^{-5}$ to approximate $\mathbf q'(\mathbf z^*, \mathbf u^*)$ at the approximate true solution solved to $10^{-16}$ accuracy.  

\autoref{fig:total_variation_e} (top) compares the error norm reduction when using ADMM, AA($m$)-ADMM and sAA($m$)-ADMM.
The convergence acceleration seen in the figure can be explained based on the spectra in
\autoref{fig:total_variation_e} (bottom).
The spectrum of $\mathbf q'(\mathbf z^*, \mathbf u^*)$ has asymptotic convergence factor
$\rho^*_{q'} = 0.976$. 
The spectrum has some complex eigenvalues.
We choose $\beta^*$ according to \autoref{prop3},
\begin{equation*}
    \beta^* = \frac{1-\sqrt{1-\rho^*_{q'}}}{1+\sqrt{1-\rho^*_{q'}}} = 0.730.
\end{equation*}

The corresponding lower bound on the optimal sAA(1)-ADMM linear convergence factor is 
\begin{equation*}
     \rho^*_{sAA(1)-ADMM} \le 1 - \sqrt{1-\rho^*_{q'}} = 0.844.
\end{equation*}
The spectral radius of the numerically computed $\Psi'$ using $\beta^*$ is given by
\begin{equation*}
\rho_{sAA(1)-ADMM}(\beta^*) = \rho(\mathbf \Psi'(\beta^*)) = 0.844 < (\rho^*_{q'})^2,
\end{equation*}
which is numerically equal to the lower bound.
It is interesting to note that it was observed numerically in \cite{desterck2020} that, for the case of sAA(1) acceleration of Alternating Least Squares for canonical tensor decomposition, for which $\mathbf q'(\mathbf x^*)$ has a complex spectrum, the lower bound in \autoref{prop3} is always achieved.

Finally, the approximately optimal $\beta_1^*$ and $\beta_2^*$ for sAA(2) are
\begin{equation*}
    \beta_1^* = 0.95, \; \beta_2^* = -0.10,
\end{equation*}
with sAA(2) linear convergence factor
\begin{equation*}
    \rho^*_{sAA(2)-ADMM} = \rho(\mathbf\Psi_2') = 0.832.
\end{equation*}
%Investigating this for sAA(1) applied to ADMM is a topic of further interest.

%\begin{figure}
%\centering
%\includegraphics[scale=0.65]{figures/total_variation_r.png}
%\caption{Total variation: comparison of error reduction using ADMM, AA($m$)-ADMM and sAA(1)-ADMM.}
%\label{fig:total_variation_r}
%\end{figure}

%-------------------------------------------------------
\subsection{Lasso problem (see, e.g., \cite{boyd2011distributed}; nonlinear and nonsmooth problem, complex spectrum)}\label{sec:lasso}
%-------------------------------------------------------
%,,,,,,,,,,,,,,,,,,,,,,,,,,,,,,,,,,,,,,,,,,,,,,,,,,,,,,,,,,,,,,,,,,,,,,,,,,,,,,
\subsubsection{Problem description}
%,,,,,,,,,,,,,,,,,,,,,,,,,,,,,,,,,,,,,,,,,,,,,,,,,,,,,,,,,,,,,,,,,,,,,,,,,,,,,,
$l_1$-regularized linear regression is also called the lasso problem:
\begin{align*}
    \min_{\mathbf x}\; & \frac{1}{2}||\mathbf A \mathbf x - \mathbf  b||_2^2+\lambda||\mathbf x||_1,
\end{align*}
where $\mathbf A\in\mathbb R^{m\times n}$ and $\mathbf b\in\mathbb R^m$ are given data, $\lambda > 0$ is a scalar regularization parameter, and $\mathbf x\in\mathbb R^n$ is the optimization variable. In typical applications, there are many more features than training examples, and the goal is to find a parsimonious model for the data \cite{boyd2011distributed}.

To apply ADMM, we solve the following constrained problem
\begin{align*}
    \min_{\mathbf x,\mathbf z}\; & \frac{1}{2}||\mathbf A \mathbf x - \mathbf  b||_2^2+\lambda||\mathbf z||_1,\\
    \text{s.t.}\; & \mathbf x - \mathbf z = 0.
\end{align*}
The scaled augmented Lagrangian is
\begin{equation*}
    L_\rho(\mathbf x, \mathbf z, \mathbf u) = \frac{1}{2}||\mathbf A \mathbf x - \mathbf  b||_2^2+\lambda||\mathbf z||_1 + \frac{\rho}{2}||\mathbf x - \mathbf z + \mathbf u||_2^2-\frac{\rho}{2}||\mathbf u||_2^2.
\end{equation*}
Therefore, we get the ADMM steps 
\begin{equation*}
\begin{cases}
    \mathbf x_{k+1} = \text{argmin}_{\mathbf x} \frac{1}{2}||\mathbf A \mathbf x - \mathbf  b||_2^2+\frac{\rho}{2}||\mathbf x - \mathbf z _k+\mathbf u_k||_2^2\\
    \mathbf z_{k+1} = \text{argmin}_{\mathbf z} \lambda||\mathbf z||_1 + \frac{\rho}{2}||\mathbf x_{k+1}-\mathbf z+\mathbf u_k||_2^2 \\
    \mathbf u_{k+1} = \mathbf u_k + \mathbf x_{k+1} - \mathbf z_{k+1},
\end{cases}
\end{equation*}
which gives
\begin{equation*}
\begin{cases}
    \mathbf x_{k+1} = (\mathbf A^T\mathbf A+\rho \mathbf I)^{-1}\left(\mathbf A^T\mathbf b + \rho(\mathbf z_k - \mathbf u_k) \right)\\
    \mathbf z_{k+1} = \text{prox}_{\frac{\lambda}{\rho}||\cdot||_1}(\mathbf x_{k+1}+\mathbf u_k), \\
    \mathbf u_{k+1} = \mathbf u_k + \mathbf x_{k+1} - \mathbf z_{k+1},
\end{cases}
\end{equation*}
where $\mathbf x_{k+1}$ can be solved efficiently as a least squares problem like in ridge regression. Since the update of $\mathbf z_{k+1}$ is nonsmooth, $\mathbf u_{k+1}$ cannot be expressed explicitly as a function of $\mathbf z_{k+1}$, and we will treat one ADMM iteration as a FPI about both variables $\mathbf z$ and $\mathbf u$ in order to apply Anderson acceleration.

\begin{figure}
\centering
\includegraphics[scale=0.65]{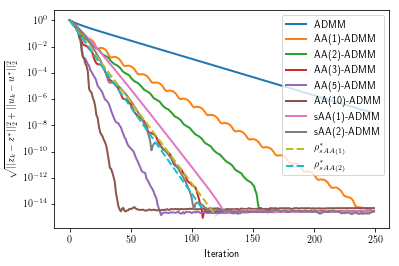}
\includegraphics[scale=0.65]{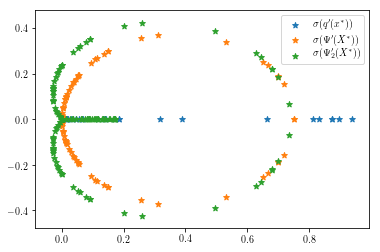}
\caption{Lasso problem (density = 0.001). (top) comparison of error reduction using ADMM, sAA($m$)-ADMM and AA($m$)-ADMM.
(bottom) Spectrum of $\mathbf q'$ of ADMM, $\mathbf\Psi'$ of sAA(1)-ADMM, and $\mathbf\Psi_2'$ of sAA(2)-ADMM.}
\label{fig:lasso_e}
\end{figure}

%\begin{figure}
%\centering
%\includegraphics[scale=0.65]{figures/eigs_lasso.png}
%\caption{Lasso problem: spectrum of ADMM iteration matrix $\mathbf q'$ and sAA($m$)-ADMM iteration matrix $\mathbf\Psi'$.}
%\label{fig:eigs_M_lasso}
%\end{figure}

%,,,,,,,,,,,,,,,,,,,,,,,,,,,,,,,,,,,,,,,,,,,,,,,,,,,,,,,,,,,,,,,,,,,,,,,,,,,,,,
\subsubsection{Parameters for the test problem}
%,,,,,,,,,,,,,,,,,,,,,,,,,,,,,,,,,,,,,,,,,,,,,,,,,,,,,,,,,,,,,,,,,,,,,,,,,,,,,,
We test our algorithms on a randomly generate sparse matrix of size $m\times n = 150\times300$ with density 0.001 and 0.01 respectively, sampled from the uniform distribution on [0,1). The $\mathbf b$ vector is sampled from the standard normal distribution. The regularization parameter $\lambda = 1$, and we pick the penalty parameter $\rho = 10$.

%,,,,,,,,,,,,,,,,,,,,,,,,,,,,,,,,,,,,,,,,,,,,,,,,,,,,,,,,,,,,,,,,,,,,,,,,,,,,,,
\subsubsection{Convergence results}
%,,,,,,,,,,,,,,,,,,,,,,,,,,,,,,,,,,,,,,,,,,,,,,,,,,,,,,,,,,,,,,,,,,,,,,,,,,,,,,
Since now the FPI is about variables $\mathbf z$ and $\mathbf u$, we will accelerate the stacked variable $[\mathbf z; \mathbf u]$. The error norm during the iteration is evaluated as
\begin{equation*}
    \mathbf e_k  = \sqrt{||\mathbf z_{k} - \mathbf z^*||_2^2 + ||\mathbf u_k-\mathbf u^*||_2^2}.
\end{equation*}
We use the first-order finite difference method with step size $h = 0.001$ to approximate $\mathbf q'(\mathbf z^*, \mathbf u^*)$ at the approximate true solution solved to $10^{-16}$ accuracy.  

\autoref{fig:lasso_e} (top) compares the error norm reduction when using ADMM, AA($m$)-ADMM and sAA($m$)-ADMM for the case when the data matrix density is 0.001.
The convergence acceleration seen in the figure can be explained based on the spectra in
\autoref{fig:lasso_e} (bottom).
The spectrum of $\mathbf q'(\mathbf z^*)$ has asymptotic convergence factor $\rho^*_{q'} = 0.938$. 
We can choose the optimal $\beta^*$ the same way as in the ridge regression problem,
\begin{equation*}
    \beta^* = \frac{1-\sqrt{1-\rho^*_{q'}}}{1+\sqrt{1-\rho^*_{q'}}} = 0.601.
\end{equation*}
The corresponding optimal sAA(1)-ADMM linear convergence factor is 
\begin{equation*}
    \rho^*_{sAA(1)-ADMM} = \rho(\mathbf\Psi') = 1 - \sqrt{1-\rho^*_{q'}} = 0.751 < (\rho^*_{q'})^2.
\end{equation*}
The approximately optimal $\beta_1^*$ and $\beta_2^*$ for sAA(2) are
\begin{equation*}
    \beta_1^* = 0.85, \; \beta_2^* = -0.10,
\end{equation*}
with sAA(2) linear convergence factor
\begin{equation*}
    \rho^*_{sAA(2)-ADMM} = \rho(\mathbf\Psi_2') = 0.737.
\end{equation*}

\begin{figure}
\centering
\includegraphics[scale=0.65]{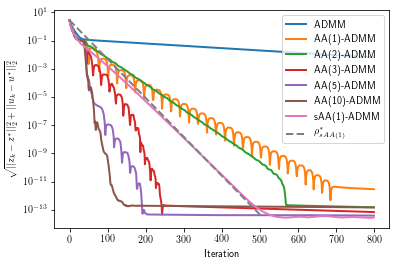}
\includegraphics[scale=0.65]{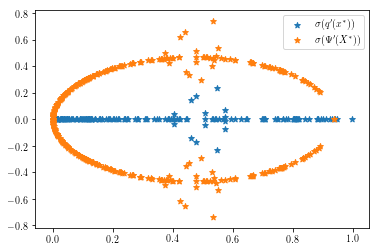}
\caption{Lasso problem (density = 0.01). (top) Comparison of error reduction using ADMM, sAA($m$)-ADMM and AA($m$)-ADMM.
(bottom) Spectrum of ADMM iteration matrix $\mathbf q'$ and sAA($1$)-ADMM iteration matrix $\mathbf\Psi'$.}
\label{fig:lasso_e_2}
\end{figure}

%\begin{figure}
%\centering
%\caption{Lasso problem (density = 0.01): }
%\label{fig:eigs_M_lasso_2}
%\end{figure}

\subsubsection{\texorpdfstring{$\mathbf q'(\mathbf x)$}{qprime} with complex eigenvalues}\label{appendix:complex_eigs}

Note that in the lasso test of \autoref{fig:lasso_e}, the eigenvalues of $\mathbf q'(\mathbf x^*)$ happen to be all real. However, this is not the case if we increase the sparsity density of data matrix $A$. For example, for a density of 0.01, $\mathbf q'(\mathbf z^*)$ has a few complex eigenvalues as shown in \autoref{fig:lasso_e_2} (bottom), where $\rho^*_{q'}=0.996$. For this case, numerical results comparing the convergence of different algorithms are shown in \autoref{fig:lasso_e_2} (top). The value of $\beta^*$ we use is
chosen according to \autoref{prop3},
\begin{equation*}
    \beta^* = \frac{1-\sqrt{1-\rho^*_{q'}}}{1+\sqrt{1-\rho^*_{q'}}} = 0.884.
\end{equation*}
The corresponding sAA(1)-ADMM linear convergence factor is 
\begin{equation*}
    \rho^*_{sAA(1)-ADMM} = \rho(\mathbf\Psi') = 1 - \sqrt{1-\rho^*_{q'}} = 0.938 < (\rho^*_{q'})^2.
\end{equation*}

\begin{figure}
\centering
\includegraphics[scale=0.65]{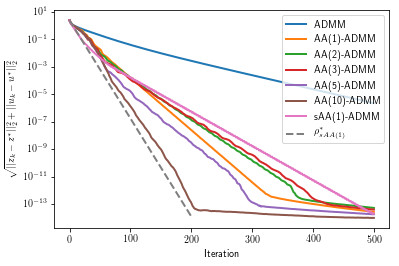}
\includegraphics[scale=0.65]{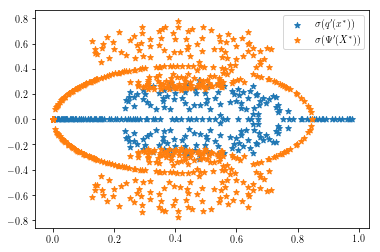}
\caption{Lasso problem (density = 0.06). (top) Comparison of error reduction using ADMM, sAA($m$)-ADMM and AA($m$)-ADMM.
(bottom) Spectrum of ADMM iteration matrix $\mathbf q'$ and sAA($1$)-ADMM iteration matrix $\mathbf\Psi'$.}
\label{fig:lasso_e_3}
\end{figure}

It is interesting to consider the situation when $\mathbf q'(\mathbf x^*)$ has complex eigenvalues with large
imaginary part. Let $\mu_+$ be the largest nonnegative real eigenvalue of $\mathbf q'(\mathbf x^*)$.
It is easy to show that, if the equality $\rho^*_{sAA(1)} = 1-\sqrt{1-\rho^*_{q'}}$ holds in \autoref{prop3}
with $\rho^*_{q'}=\mu_+$ and $\beta^* = \frac{1-\sqrt{1-\rho^*_{q'}}}{1+\sqrt{1-\rho^*_{q'}}}$,
then the rightmost point of the circle of \autoref{prop:circle} is the image of
$\mu_+$ under the mapping from $\mu$ to $\lambda$ defined by (\ref{eqn:sAA2_lambda_mu_relation}),
and this point determines $\rho(\mathbf \Psi')=\rho^*_{sAA(1)}$.
According to Corollary S.1 in the supplementary materials of \cite{desterck2020}, this also holds
when $\rho^*_{q'}>\mu_+$ and $\rho^*_{sAA(1)} = 1-\sqrt{1-\mu_+}$, with
$\beta^* = \frac{1-\sqrt{1-\mu_+}}{1+\sqrt{1-\mu_+}}$. In these cases, the spectral radius
of $\mathbf \Psi'$ for sAA(1) with optimal weight is determined by the mapped eigenvalue
of $\mu_+$, which is the rightmost point of the circle, and the complex eigenvalues of
$\mathbf q'(\mathbf x^*)$ do not influence $\rho(\mathbf \Psi')$.
However, when $\mathbf q'(\mathbf x^*)$ has complex eigenvalues with large
imaginary part, these eigenvalues may be mapped to eigenvalues $\lambda$
of $\mathbf \Psi'$ that are sufficiently far outside the circle of \autoref{prop:circle}
to determine the spectral radius of $\mathbf \Psi'$.
In this case, we cannot determine the optimal $\beta^*$ and $\rho^*_{sAA(1)}$
by the expressions (with equality) in \autoref{prop3} or Corollary S.1
in the supplementary materials of \cite{desterck2020}.
We now give an example demonstrating this.
We consider the lasso example with density = 0.06.
\autoref{fig:lasso_e_3} (bottom) plots the distribution of eigenvalues for both $\mathbf q'(\mathbf z^*, \mathbf u^*)$ and $\mathbf\Psi'$, where we have used $\beta = \frac{1-\sqrt{1-\rho^*_{q'}}}{1+\sqrt{1-\rho^*_{q'}}}$ in sAA(1)-ADMM. We can see that the largest eigenvalues of $\mathbf\Psi'$ induced by complex eigenvalues of $\mathbf q'$ (those that are not lying on the circle) have a larger modulus (=0.944) than the largest-size eigenvalue induced by the real eigenvalues of $\mathbf q'$, which is of size $1 - \sqrt{1-\rho^*_{q'}} = 0.848$ (since for this example it still holds that $\rho^*_{q'}=\mu_+$).
Hence, complex eigenvalues of $\mathbf q'$ dominate the spectrum of $\mathbf \Psi'$ and the
equality in \autoref{prop3} does not hold, since it requires that the largest-size eigenvalue of $\mathbf\Psi'$ comes from real eigenvalues of $\mathbf q'$. 
%To see this, remember that we use a quadratic equation and a real eigenvalue of $\mathbf q'$ to get two eigenvalues of $\mathbf\Psi'$. 
This observation matches with the numerical results shown in \autoref{fig:lasso_e_3} (top), where the convergence of the sAA(1) algorithm using $\beta = (1 - \sqrt{1-\rho^*_{q'}})/(1 + \sqrt{1-\rho^*_{q'}})$ from \autoref{prop3}
does not match the convergence factor $1 - \sqrt{1-\rho^*_{q'}}$ that would correspond to \autoref{prop3} if
equality were to hold. (Note that for the previous test with density 0.01 we do get a close match
(see \autoref{fig:lasso_e_2} (top)).)
We note that in a case like the one from \autoref{fig:lasso_e_3}, sAA(1) may generate divergent results when using
$\beta = (1 - \sqrt{1-\rho^*_{q'}})/(1 + \sqrt{1-\rho^*_{q'}})$ since this is not the correct optimal $\beta^*$.
Finding the optimal coefficient $\beta^*$ for sAA(1) in this scenario is an open question that needs more investigation.

%-------------------------------------------------------
\subsection{Nonnegative least squares (see, e.g., \cite{fu2019anderson}; nonlinear problem with inequality constraint)}
%-------------------------------------------------------
%,,,,,,,,,,,,,,,,,,,,,,,,,,,,,,,,,,,,,,,,,,,,,,,,,,,,,,,,,,,,,,,,,,,,,,,,,,,,,,
\subsubsection{Problem description}
%,,,,,,,,,,,,,,,,,,,,,,,,,,,,,,,,,,,,,,,,,,,,,,,,,,,,,,,,,,,,,,,,,,,,,,,,,,,,,,
The nonnegative least squares problem is 
\begin{equation*}
    \min_{\mathbf x} ||\mathbf F\mathbf x - \mathbf g||_2^2, \quad \text{s.t.} \quad \mathbf x \geq 0,
\end{equation*}
where $\mathbf x \in \mathbb R^n$ is the variable, and $\mathbf F \in \mathbb R^{m\times n}$ and $\mathbf g\in\mathbb R^m$ are problem data. We can integrate the nonnegativity constraint into the objective function and rewrite the problem as
\begin{align*}
    \min_{\mathbf x, \mathbf z}\; &||\mathbf F\mathbf x - \mathbf g||_2^2 + \mathcal I_{\mathbb R_+^n}(\mathbf z),\\
    \text{s.t.}\; &\mathbf x - \mathbf z = 0,
\end{align*}
where $\mathcal{I}_{\mathbb R_+^n}$ is the indicator function defined as
\begin{equation*}
    \mathcal{I}_{\mathbb R_+^{n}}(\mathbf z) =
    \begin{cases}
        0, & \mathbf z \geq 0\\
        +\infty, & \text{otherwise}.
    \end{cases}
\end{equation*}

The scaled augmented Lagrangian of this problem is
\begin{equation*}
    L_{\rho}(\mathbf x, \mathbf z, \mathbf u) = ||\mathbf F \mathbf x - \mathbf g||_2^2 + \mathcal{I}_{\mathbb R_+^{n}}(\mathbf z) + \frac{\rho}{2}||\mathbf x - \mathbf z + \mathbf u||_2^2 - \frac{\rho}{2}||\mathbf u||_2^2.
\end{equation*}
The ADMM steps on this problem are:
\begin{equation*}
\begin{cases}
    \mathbf x_{k+1} = \text{argmin}_{\mathbf x} ||\mathbf F\mathbf x - \mathbf g||_2^2 + \frac{\rho}{2} ||\mathbf x - \mathbf z_k + \mathbf u_k||_2^2\\
    \mathbf z_{k+1} = \text{argmin}_{\mathbf z} \mathcal{I}_{\mathbb R_+^{n}}(\mathbf z) + \frac{\rho}{2} ||\mathbf x_{k+1}  + \mathbf u_k - \mathbf z||_2^2 \\
    \mathbf u_{k+1} = \mathbf u_k + \mathbf x_{k+1} - \mathbf z_{k+1}
\end{cases}
\end{equation*}
where the first step for $\mathbf x_{k+1}$ is the proximal operator of the $l_2$-norm.
% and can be solved from the least square problem
%\begin{equation*}
%    \mathbf x_{k+1} = \text{argmin}_{\mathbf x}\left|\left| 
%    \begin{bmatrix}
%    F \\
%    \sqrt{\frac{\rho}{2}}I
%    \end{bmatrix} x -
%    \begin{bmatrix}
%    g \\
%    \sqrt{\frac{\rho}{2}}(z_k - u_k)
%    \end{bmatrix}
%    \right|\right|_2^2.
%\end{equation*}
The second step is just the proximal operator of the indicator function, which is equivalent to the projection operator
\begin{equation*}
    \mathbf z_{k+1} = \frac{1}{\rho}\Uppi_{\mathbb R_+^n}(\mathbf x_{k+1} + u_k).
\end{equation*}

\begin{figure}
\centering
\includegraphics[scale=0.65]{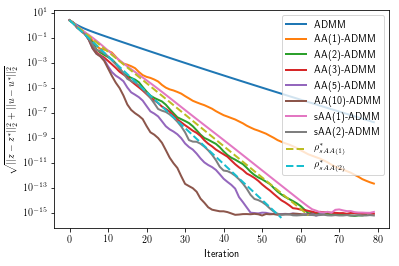}
\includegraphics[scale=0.65]{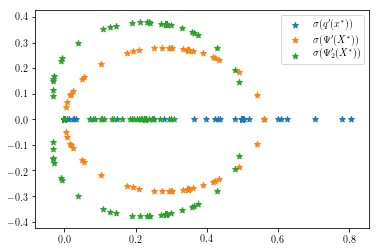}
\caption{Nonnegative least squares. (top) Comparison of error reduction using ADMM, sAA($m$)-ADMM and AA($m$)-ADMM.
(bottom) Spectrum of $\mathbf q'$ of ADMM, $\mathbf\Psi'$ of sAA(1)-ADMM, and $\mathbf\Psi_2'$ of sAA(2)-ADMM.}
\label{fig:nnls_e}
\end{figure}

%\begin{figure}
%\centering
%\caption{Nonnegative least squares: }
%\label{fig:eigs_nnls}
%\end{figure}

%,,,,,,,,,,,,,,,,,,,,,,,,,,,,,,,,,,,,,,,,,,,,,,,,,,,,,,,,,,,,,,,,,,,,,,,,,,,,,,
\subsubsection{Parameters for the test problem}
%,,,,,,,,,,,,,,,,,,,,,,,,,,,,,,,,,,,,,,,,,,,,,,,,,,,,,,,,,,,,,,,,,,,,,,,,,,,,,,
We test our algorithms on a randomly generated sparse matrix of size $m\times n = 150\times300$ with density 0.001, sampled from the standard normal distribution. The $g$ vector is sampled from the standard normal distribution. The augmented Lagrangian penalty parameter is $\rho = 2$.

%,,,,,,,,,,,,,,,,,,,,,,,,,,,,,,,,,,,,,,,,,,,,,,,,,,,,,,,,,,,,,,,,,,,,,,,,,,,,,,
\subsubsection{Convergence results}
%,,,,,,,,,,,,,,,,,,,,,,,,,,,,,,,,,,,,,,,,,,,,,,,,,,,,,,,,,,,,,,,,,,,,,,,,,,,,,,
We use the first-order finite difference method with step size $h = 0.001$ to approximate $\mathbf q'(\mathbf z^*, \mathbf u^*)$ at the approximate true solution solved to $10^{-16}$ accuracy.  

\autoref{fig:nnls_e} (top) compares the error norm reduction when using ADMM, AA($m$)-ADMM and sAA($m$)-ADMM.
The convergence acceleration seen in the figure can be explained based on the spectra in
\autoref{fig:nnls_e} (bottom).
The spectrum of $\mathbf q'(\mathbf z^*, \mathbf u^*)$ has asymptotic convergence factor
$\rho^*_{q'} = 0.806$. We can choose the optimal $\beta^*$ the same way as in the ridge regression problem 
\begin{equation*}
    \beta^* = \frac{1-\sqrt{1-\rho^*_{q'}}}{1+\sqrt{1-\rho^*_{q'}}} = 0.389.
\end{equation*}
The corresponding optimal sAA(1)-ADMM linear convergence factor is 
\begin{equation*}
    \rho^*_{sAA(1)-ADMM} = \rho(\mathbf\Psi') = 1 - \sqrt{1-\rho^*_{q'}} = 0.560 < (\rho^*_{q'})^2.
\end{equation*}
The approximately optimal $\beta_1^*$ and $\beta_2^*$ for sAA(2) are
\begin{equation*}
    \beta_1^* = 0.70, \; \beta_2^* = -0.10,
\end{equation*}
with sAA(2) linear convergence factor
\begin{equation*}
    \rho^*_{sAA(2)-ADMM} = \rho(\mathbf\Psi_2') = 0.516.
\end{equation*}
\subsection{Constrained logistic regression (see, e.g., \cite{mai2019anderson}; nonlinear problem with box constraint)}
%-------------------------------------------------------
%,,,,,,,,,,,,,,,,,,,,,,,,,,,,,,,,,,,,,,,,,,,,,,,,,,,,,,,,,,,,,,,,,,,,,,,,,,,,,,
\subsubsection{Problem description}
%,,,,,,,,,,,,,,,,,,,,,,,,,,,,,,,,,,,,,,,,,,,,,,,,,,,,,,,,,,,,,,,,,,,,,,,,,,,,,,
The constrained regularized logistic regression adds a constraint on $||\mathbf x||_{\infty}$ to the regularized logistic regression problem that we have already discussed:
\begin{align*}
    \min_{\mathbf x}\; &\frac{1}{m}\sum_{i=1}^m\log(1 + \exp(-y_i(\mathbf a_i^T\mathbf w+\mathbf c))) + \lambda ||\mathbf x||_2^2,\\
    \text{s.t.}\;& ||\mathbf x||_{\infty} \leq 1.
\end{align*}
To apply ADMM, we rewrite this problem as
\begin{align*}
    \min_{\mathbf x, \mathbf z}\; &\frac{1}{m}\sum_{i=1}^m\log(1 + \exp(-y_i(a_i^T\mathbf w+\mathbf c))) + \lambda ||\mathbf x||_2^2 + \mathcal I_{\Omega}(\mathbf z),\\
    \text{s.t.}\;& \mathbf x -\mathbf z = 0,
\end{align*}
where $\Upomega = \{\mathbf x : ||\mathbf x||_{\infty}\leq 1\}$. This gives the augmented Lagrangian
\begin{equation*}
    L_\rho(\mathbf x, \mathbf z, \mathbf u) = \frac{1}{m}\sum_{i=1}^m\log(1 + \exp(-y_i(\mathbf a_i^T\mathbf w+\mathbf c))) + \lambda ||\mathbf x||_2^2 + \mathcal I_{\Omega}(\mathbf z) + \frac{\rho}{2}||\mathbf x - \mathbf z + \mathbf u||_2^2 - \frac{\rho}{2}||\mathbf u||_2^2.
\end{equation*}
Hence, we get the ADMM steps
\begin{equation*}
\begin{cases}
    \mathbf x_{k+1} = \text{argmin}_{\mathbf x} \frac{1}{m}\sum_{i=1}^m\log(1 + \exp(-y_i(\mathbf a_i^T\mathbf w+\mathbf c))) + \lambda ||\mathbf x||_2^2 +  \frac{\rho}{2}||\mathbf x - \mathbf z _k+\mathbf u_k||_2^2\\
    \mathbf z_{k+1} = \text{argmin}_{\mathbf z} \mathcal I_{\Omega}(\mathbf z) + \frac{\rho}{2}||\mathbf x_{k+1}-\mathbf z+\mathbf u_k||_2^2 \\
    \mathbf u_{k+1} = \mathbf u_k + \mathbf x_{k+1} - \mathbf z_{k+1}.
\end{cases}
\end{equation*}
Like before, we use Newton's method to solve for $\mathbf x_{k+1}$. For $\mathbf z_{k+1}$, since the proximal operation of an indicator function is just a projection, we have
\begin{equation*}
    \mathbf z_{k+1} = \frac{1}{\rho}\Uppi_{\Upomega}(\mathbf x_{k+1}+\mathbf u_k),
\end{equation*}
which is
\begin{equation*}
    [\mathbf z_{k+1}]_j = 
    \begin{cases}
        \frac{1}{\rho}, & [\mathbf x_{k+1}+\mathbf u_k]_j \in [1, \infty)\\
        \frac{1}{\rho}[\mathbf x_{k+1}+\mathbf u_k]_j, & [\mathbf x_{k+1}+\mathbf u_k]_j \in (-1,1)\\
        -\frac{1}{\rho}, & [\mathbf x_{k+1}+\mathbf u_k]_j \in (-\infty, -1].
    \end{cases}
\end{equation*}

\begin{figure}
\centering
\includegraphics[scale=0.65]{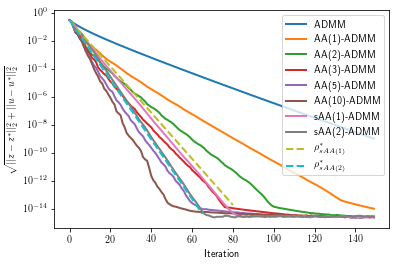}
\includegraphics[scale=0.65]{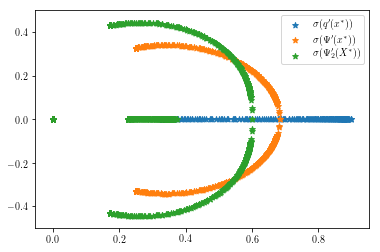}
\caption{Constrained regularized logistic regression. (top) Comparison of error reduction using ADMM, sAA($m$)-ADMM and AA($m$)-ADMM.
(bottom) Spectrum of ADMM iteration matrix $\mathbf q'_*$ and sAA($1$)-ADMM iteration matrix $\mathbf \Psi'_*$.}
\label{fig:constr_reglog_e}
\end{figure}

%\begin{figure}
%\centering
%\caption{Constrained regularized logistic regression: }
%\label{fig:eigs_constr_reglog}
%\end{figure}

%,,,,,,,,,,,,,,,,,,,,,,,,,,,,,,,,,,,,,,,,,,,,,,,,,,,,,,,,,,,,,,,,,,,,,,,,,,,,,,
\subsubsection{Parameters for the test problem}
%,,,,,,,,,,,,,,,,,,,,,,,,,,,,,,,,,,,,,,,,,,,,,,,,,,,,,,,,,,,,,,,,,,,,,,,,,,,,,,
We use the same sample data from the Madelon data set as in Section \ref{sec:reglog}. The regularization and penalty parameters are $\lambda = 2$ and $\rho = 10$
respectively, as in Section \ref{sec:reglog}.

%,,,,,,,,,,,,,,,,,,,,,,,,,,,,,,,,,,,,,,,,,,,,,,,,,,,,,,,,,,,,,,,,,,,,,,,,,,,,,,
\subsubsection{Convergence results}
%,,,,,,,,,,,,,,,,,,,,,,,,,,,,,,,,,,,,,,,,,,,,,,,,,,,,,,,,,,,,,,,,,,,,,,,,,,,,,,
We use the first-order finite difference method with step size $h = 0.001$ to approximate $\mathbf q'(\mathbf z^*, \mathbf u^*)$ at the approximate true solution solved to $10^{-16}$ accuracy.  

\autoref{fig:constr_reglog_e} (top) compares the error norm reduction when using ADMM, AA($m$)-ADMM and sAA($m$)-ADMM.
The convergence acceleration seen in the figure can be explained based on the spectra in
\autoref{fig:constr_reglog_e} (bottom).
The spectrum of $\mathbf q'(\mathbf z^*, \mathbf u^*)$ has asymptotic convergence factor
$\rho^*_{q'} = 0.900$. We can choose the optimal $\beta^*$ the same way as in the ridge regression problem 
\begin{equation*}
    \beta^* = \frac{1-\sqrt{1-\rho^*_{q'}}}{1+\sqrt{1-\rho^*_{q'}}} = 0.519.
\end{equation*}
The corresponding optimal sAA(1)-ADMM linear convergence factor is 
\begin{equation*}
    \rho^*_{sAA(1)-ADMM} = \rho(\mathbf\Psi') = 1 - \sqrt{1-\rho^*_{q'}} = 0.684 < (\rho^*_{q'})^2.
\end{equation*}
The approximately optimal $\beta_1^*$ and $\beta_2^*$ for sAA(2) are
\begin{equation*}
    \beta_1^* = 0.90, \; \beta_2^* = -0.15,
\end{equation*}
with sAA(2) linear convergence factor
\begin{equation*}
    \rho^*_{sAA(2)-ADMM} = \rho(\mathbf\Psi_2') = 0.612.
\end{equation*}

%\begin{figure}
%\centering
%\includegraphics[scale=0.65]{constr_reglog_r.png}
%\caption{Constrained regularized logistic regression: comparison of residual reduction using ADMM, sAA(1)-ADMM and AA($m$)-ADMM}
%\label{fig:constr_reglog_r}
%\end{figure}

%===================================
\section{Conclusions}
\label{sec:conc}
%===================================
This paper has discussed a strategy for computing the optimal asymptotic convergence factor
of stationary Anderson acceleration applied to ADMM, for the case where ADMM by itself converges linearly.
Based on the spectra of $\mathbf q'(\mathbf x^*))$ and sAA($m$)-ADMM we have provided new insight
into how the acceleration is achieved.
This approach, based on theoretical results from \cite{desterck2020}, finds numerically
that convergence factors of the stationary form of Anderson acceleration with coefficients that are chosen to
make the convergence factors optimal, provide a useful prediction for the asymptotic convergence speed of 
non-stationary AA with finite window size, which is the method used in practice.
As discussed in \cite{desterck2020}, this is intuitively reasonable: the nonstationary AA does not use these
\emph{globally optimal} stationary coefficients, but rather performs a \emph{local optimization} of the coefficients
in every step $k$ by solving least squares problem (\ref{eqn:anderson_coeff_k-i_k-i-1}).
As $\mathbf x$ approaches $\mathbf x^*$ in the asymptotic regime and $\mathbf q'(\mathbf x)$ approaches
$\mathbf q'(\mathbf x^*)$, it is not unreasonable to expect the convergence behavior of AA with locally-optimal
$\beta_i^{(k)}$ weights to be similar to the behavior of sAA with
weights that are, based on $\mathbf q'(\mathbf x^*)$, globally optimal in obtaining the best asymptotic
convergence rate.
This is indeed what we have observed numerically in this paper for AA applied to ADMM.

The case of sAA with $m=1$ is easy to analyze and directly leads to the simple
analytical prediction formulas of \autoref{prop2} and \autoref{prop3} for the optimal convergence
factors $\rho^*_{sAA(1)}$, see \cite{desterck2020}. 
While our numerical results show that $\rho^*_{sAA(1)-ADMM}$ is a useful prediction
for $\rho_{AA(m)-ADMM}$ also when $m>1$, it is clear that computing
$\rho^*_{sAA(m)-ADMM}$ for $m>1$ is also of interest.
As we have illustrated, for $m\ge2$ the optimal $\rho^*_{sAA(m)-ADMM}$
can be obtained by optimization \cite{desterck2020}, but the lack of analytical results is
an interesting avenue for further research, for example, on how the optimal 
$\rho^*_{sAA(m)-ADMM}$ depends on $m$.

The similarity in asymptotic convergence behavior between AA and optimal sAA
allows us to understand the acceleration power of AA in terms of how it reshapes
convergence spectra in our numerical tests, in ways that are very similar to how
GMRES for linear systems accelerates convergence depending on the spectral
and eigenspace properties of the GMRES preconditioner (see \cite{desterck2020}
for a detailed discussion of this analogy).
%; indeed, ADMM can be viewed as an effective
%nonlinear preconditioner for AA (see \cite{desterck2020} for a detailed discussion of
%this analogy).

The similarity between AA and optimal sAA convergence factors also provides
a prediction for convergence acceleration by AA, which is especially
useful since the quest for linear asymptotic convergence bounds for AA
with finite window size has been elusive, due to the AA coefficients changing in every
iteration. This similarity may also inspire theoretical approaches for finding
asymptotic convergence factor bounds for AA with finite window size.

Of course, besides providing useful insight, our approach for estimating AA convergence
factors is not really practical, since $\rho(\mathbf q'(\mathbf x^*))$ needs to be known or computed 
to compute the optimal $\rho^*_{sAA(1)-ADMM}$.
However, if an upper bound for $\rho(\mathbf q'(\mathbf x^*))$ is known,
then an upper bound for the optimal sAA(1)-ADMM convergence
factor, $\rho^*_{sAA(1)-ADMM}$, can directly be obtained from the formulas
in \autoref{prop2} and \autoref{prop3}.
In preconditioned GMRES for linear systems, depending on the problem,
such upper bounds for $\rho(\mathbf q'(\mathbf x^*))$ can often be derived \cite{desterck2020}.
They may, for example, depend on problem parameters or problem sizes, and
for many linear problems GMRES preconditioners have been found that provably
lead to favorable convergence bounds independent from, or only weakly dependent
on, parameters that characterize the difficulty or conditioning of the problem.
Similarly, it may be of practical use to pursue this for various ADMM applications,
since it may lead to convergence factor bound predictions
for AA applied to ADMM with favorable dependence on problem parameters.

%\begin{acknowledgements}
%If you'd like to thank anyone, place your comments here
%and remove the percent signs.
%\end{acknowledgements}

% Authors must disclose all relationships or interests that 
% could have direct or potential influence or impart bias on 
% the work: 
%

% BibTeX users please use one of
%\bibliographystyle{spbasic}      % basic style, author-year citations
\bibliographystyle{spmpsci}      % mathematics and physical sciences
\bibliography{AA-ADMM}   % name your BibTeX data base

\begin{thebibliography}{10}
\providecommand{\url}[1]{{#1}}
\providecommand{\urlprefix}{URL }
\expandafter\ifx\csname urlstyle\endcsname\relax
  \providecommand{\doi}[1]{DOI~\discretionary{}{}{}#1}\else
  \providecommand{\doi}{DOI~\discretionary{}{}{}\begingroup
  \urlstyle{rm}\Url}\fi

\bibitem{anderson1965iterative}
Anderson, D.G.: Iterative procedures for nonlinear integral equations.
\newblock Journal of the ACM (JACM) \textbf{12}(4), 547--560 (1965)

\bibitem{boley2012linear}
Boley, D.: Linear convergence of {ADMM} on a model problem.
\newblock Department of Computer Science and Engineering, University of
  Minnesota, TR pp. 12--009 (2012)

\bibitem{boyd2011distributed}
Boyd, S., Parikh, N., Chu, E.: Distributed optimization and statistical
  learning via the alternating direction method of multipliers.
\newblock Now Publishers Inc (2011)

\bibitem{davis2017faster}
Davis, D., Yin, W.: Faster convergence rates of relaxed {Peaceman-Rachford} and
  {ADMM} under regularity assumptions.
\newblock Mathematics of Operations Research \textbf{42}(3), 783--805 (2017)

\bibitem{sterck2012nonlinear}
{De Sterck}, H.: A nonlinear {GMRES} optimization algorithm for canonical
  tensor decomposition.
\newblock SIAM J. Scientific Computing \textbf{34}(3), A1351--A1379 (2012)

\bibitem{desterck2020}
De~Sterck, H., He, Y.: On the asymptotic linear convergence speed of {Anderson}
  acceleration, nestrov acceleration and nonlinear {GMRES}.
\newblock to appear, SIAM Journal on Scientific Computing;
  \url{https://arxiv.org/abs/2007.01996}  (2020)

\bibitem{deng2016global}
Deng, W., Yin, W.: On the global and linear convergence of the generalized
  alternating direction method of multipliers.
\newblock Journal of Scientific Computing \textbf{66}(3), 889--916 (2016)

\bibitem{franca2018admm}
Franca, G., Robinson, D.P., Vidal, R.: {ADMM} and accelerated {ADMM} as
  continuous dynamical systems.
\newblock \url{https://arxiv.org/abs/1805.06579}  (2018)

\bibitem{francca2018dynamical}
Fran{\c{c}}a, G., Robinson, D.P., Vidal, R.: A dynamical systems perspective on
  nonsmooth constrained optimization.
\newblock \url{https://arxiv.org/abs/1808.04048}  (2018)

\bibitem{fu2019anderson}
Fu, A., Zhang, J., Boyd, S.: Anderson accelerated {Douglas-Rachford} splitting.
\newblock \url{https://arxiv.org/abs/1908.11482}  (2019)

\bibitem{ghadimi2014optimal}
Ghadimi, E., Teixeira, A., Shames, I., Johansson, M.: Optimal parameter
  selection for the alternating direction method of multipliers ({ADMM}):
  quadratic problems.
\newblock IEEE Transactions on Automatic Control \textbf{60}(3), 644--658
  (2014)

\bibitem{goldstein2014fast}
Goldstein, T., O'Donoghue, B., Setzer, S., Baraniuk, R.: Fast alternating
  direction optimization methods.
\newblock SIAM Journal on Imaging Sciences \textbf{7}(3), 1588--1623 (2014)

\bibitem{he20121}
He, B., Yuan, X.: On the $\mathcal{O}(1/n)$ convergence rate of the
  {Douglas--Rachford} alternating direction method.
\newblock SIAM Journal on Numerical Analysis \textbf{50}(2), 700--709 (2012)

\bibitem{he2015non}
He, B., Yuan, X.: On non-ergodic convergence rate of {Douglas--Rachford}
  alternating direction method of multipliers.
\newblock Numerische Mathematik \textbf{130}(3), 567--577 (2015)

\bibitem{hong2017linear}
Hong, M., Luo, Z.Q.: On the linear convergence of the alternating direction
  method of multipliers.
\newblock Mathematical Programming \textbf{162}(1-2), 165--199 (2017)

\bibitem{kadkhodaie2015accelerated}
Kadkhodaie, M., Christakopoulou, K., Sanjabi, M., Banerjee, A.: Accelerated
  alternating direction method of multipliers.
\newblock In: Proceedings of the 21th ACM SIGKDD international conference on
  knowledge discovery and data mining, pp. 497--506 (2015)

\bibitem{lions1979splitting}
Lions, P.L., Mercier, B.: Splitting algorithms for the sum of two nonlinear
  operators.
\newblock SIAM Journal on Numerical Analysis \textbf{16}(6), 964--979 (1979)

\bibitem{mai2019anderson}
Mai, V.V., Johansson, M.: Anderson acceleration of proximal gradient methods.
\newblock \url{https://arxiv.org/abs/1910.08590}  (2019)

\bibitem{mitchell2020nesterov}
Mitchell, D., Ye, N., De~Sterck, H.: {Nesterov} acceleration of alternating
  least squares for canonical tensor decomposition: {Momentum} step size
  selection and restart mechanisms.
\newblock Numerical Linear Algebra with Applications p. e2297 (2020)

\bibitem{nishihara2015general}
Nishihara, R., Lessard, L., Recht, B., Packard, A., Jordan, M.I.: A general
  analysis of the convergence of {ADMM}.
\newblock \url{https://arxiv.org/abs/1502.02009}  (2015)

\bibitem{ortega2000iterative}
Ortega, J.M., Rheinboldt, W.C.: Iterative solution of nonlinear equations in
  several variables.
\newblock SIAM (2000)

\bibitem{peng2018anderson}
Peng, Y., Deng, B., Zhang, J., Geng, F., Qin, W., Liu, L.: Anderson
  acceleration for geometry optimization and physics simulation.
\newblock ACM Transactions on Graphics (TOG) \textbf{37}(4), 1--14 (2018)

\bibitem{poon2019trajectory}
Poon, C., Liang, J.: Trajectory of alternating direction method of multipliers
  and adaptive acceleration.
\newblock In: Advances in Neural Information Processing Systems, pp. 7355--7363
  (2019)

\bibitem{toth2015convergence}
Toth, A., Kelley, C.: Convergence analysis for {Anderson} acceleration.
\newblock SIAM Journal on Numerical Analysis \textbf{53}(2), 805--819 (2015)

\bibitem{walker2011anderson}
Walker, H.F., Ni, P.: Anderson acceleration for fixed-point iterations.
\newblock SIAM Journal on Numerical Analysis \textbf{49}(4), 1715--1735 (2011)

\bibitem{zhang2019accelerating}
Zhang, J., Peng, Y., Ouyang, W., Deng, B.: Accelerating {ADMM} for efficient
  simulation and optimization.
\newblock ACM Transactions on Graphics (TOG) \textbf{38}(6), 1--21 (2019)

\bibitem{zhang2018gmres}
Zhang, R.Y., White, J.K.: {GMRES}-accelerated {ADMM} for quadratic objectives.
\newblock SIAM Journal on Optimization \textbf{28}(4), 3025--3056 (2018)

\end{thebibliography}

\begin{appendices}

%\begin{figure}
%\centering
%\caption{Ridge regression: spectrum of ADMM iteration matrix $\mathbf q'$ and sAA($m$)-ADMM iteration matrix $\mathbf\Psi'_m$.}
%\label{fig:eigs_ridge_sAA3}
%\end{figure}

%\begin{figure}
%\centering
%\caption{Regularized logistic regression: spectrum of ADMM iteration matrix $\mathbf q'$ and sAA($m$)-ADMM iteration matrix $\mathbf\Psi'_m$.}
%\label{fig:eigs_logreg_sAA3}
%\end{figure}

\section{Derivative of $\mathbf q(\mathbf x)$ for \texorpdfstring{$l_1$}{L1}-regularized problems}\label{appendix:l1-norm}
We mentioned in the main text that although the $l_1$-regularized least squares problem is nonsmooth, the FPI representation $\mathbf q(\cdot)$ of ADMM can be differentiable at the true solution $\mathbf z^*$. To see this, consider the following simple scalar example
\begin{equation*}
    \min_{x\in\mathbb R} f(\mathbf x) := \frac{1}{2}x^2 + |x|.
\end{equation*}
Clearly, the objective function is non-differentiable at $x = 0$ and the optimum is also achieved at $x = 0$. The equivalent split form is
\begin{equation*}
\begin{aligned}
    \min_{x, z}\; &\frac{1}{2}x^2 + |z|, \\
    \text{s.t. }\; & x - z  = 0.
\end{aligned}
\end{equation*}
From the ADMM update
\begin{equation*}
\begin{cases}
   x_{k+1} = \frac{\rho}{1+\rho}( z_k - u_k),\\
    z_{k+1} = \textrm{prox}_{\frac{1}{\rho}|\cdot|}(x_{k+1}+ u_k), \\
    u_{k+1} = u_k + x_{k+1}- z_{k+1},
\end{cases}
\end{equation*}
we can get the FPI representation $(z_{k+1}, u_{k+1}) = \mathbf q(z_k, u_k)$, where

\begin{equation*}
\begin{bmatrix}
    z_{k+1}\\
    u_{k+1}
\end{bmatrix} =
    \begin{cases}
        \begin{bmatrix}  % case 1
            \frac{\rho}{1+\rho} & \frac{1}{1+\rho}\\
            0                   & 0
        \end{bmatrix}
        \begin{bmatrix}
            z_k\\
            u_k
        \end{bmatrix}
        +
        \begin{bmatrix}
            -\frac{1}{\rho}\\
            \frac{1}{\rho}
        \end{bmatrix}
        & \text{if } \frac{\rho}{1+\rho}z_k + \frac{1}{1+\rho}u_k > \frac{1}{\rho},\\
        \\
        \begin{bmatrix}  % case 2
            0                   & 0\\
            \frac{\rho}{1+\rho} & \frac{1}{1+\rho}
        \end{bmatrix}
        \begin{bmatrix}
            z_k\\
            u_k
        \end{bmatrix}
        & \text{if } \left| \frac{\rho}{1+\rho}z_k + \frac{1}{1+\rho}u_k \right| \leq \frac{1}{\rho},\\
        \\
        \begin{bmatrix}  % case 3
            \frac{\rho}{1+\rho} & \frac{1}{1+\rho}\\
            0                   & 0
        \end{bmatrix}
        \begin{bmatrix}
            z_k\\
            u_k
        \end{bmatrix}
        +
        \begin{bmatrix}
            \frac{1}{\rho}\\
            -\frac{1}{\rho}
        \end{bmatrix}
        & \text{if } \frac{\rho}{1+\rho}z_k + \frac{1}{1+\rho}u_k < -\frac{1}{\rho},
    \end{cases}
\end{equation*}
which is a nonsmooth function. Hence, we have
\begin{equation*}
\mathbf q'\left(z, u\right)
=
    \begin{cases}
        \begin{bmatrix}  % case 1
            \frac{\rho}{1+\rho} & \frac{1}{1+\rho}\\
            0                   & 0
        \end{bmatrix}
        & \text{if } \frac{\rho}{1+\rho}z + \frac{1}{1+\rho}u > \frac{1}{\rho},\\
        \\
        \begin{bmatrix}  % case 2
            0                   & 0\\
            \frac{\rho}{1+\rho} & \frac{1}{1+\rho}
        \end{bmatrix}
        & \text{if } \left| \frac{\rho}{1+\rho}z + \frac{1}{1+\rho}u \right| \leq \frac{1}{\rho},\\
        \\
        \begin{bmatrix}  % case 3
            \frac{\rho}{1+\rho} & \frac{1}{1+\rho}\\
            0                   & 0
        \end{bmatrix}
        & \text{if } \frac{\rho}{1+\rho}z + \frac{1}{1+\rho}u < -\frac{1}{\rho}.\\
    \end{cases}
\end{equation*}
Because the optimal solution is $z^* = u^* = 0$, we see that $\mathbf q'(z^*, u^*)$ exists and 
\begin{equation*}
\mathbf q'\left(
    z^*,
    u^*
\right)
=
\begin{bmatrix}  % case 2
    0                   & 0\\
    \frac{\rho}{1+\rho} & \frac{1}{1+\rho}
\end{bmatrix}.
\end{equation*}
From this example, we can see that even when the objective function is nonsmooth, the FPI representation of ADMM for solving the problem can still be differentiable at the optimal solution.
In this example this is the case as long as $\frac{\rho}{1+\rho}z^* + \frac{1}{1+\rho}u^* \pm \frac{1}{\rho} \neq 0$, where $\frac{\rho}{1+\rho}z + \frac{1}{1+\rho}u \pm \frac{1}{\rho}$ is obtained from the proximal operation of the $z$-update at its nondifferentiable point $z = 0$. We see that the soft-thresholding operation spreads out the nondifferentiable point at $z = 0$ in the original problem to two lines in the $z,u$ plane. From the optimality conditions
\begin{equation*}
    x^*-z^* = 0, \quad x^* + \rho u^* = 0,
\end{equation*}
we can get
\begin{equation*}
    z^* + \rho u^* = 0.
\end{equation*}
Therefore, only when
\begin{equation*}
    z^*=\pm\frac{1}{1-\rho}, \quad u^* = \mp\frac{1}{\rho(1-\rho)},
\end{equation*}
is $\mathbf q(\cdot)$ not differentiable at the true solution, 
but the true solution is $x^*=z^*=0$ in this example.

We can generalize this observation to multi-dimensional problems. For example, 
for the total variation problem of Section \ref{sec:totvar}, we get
\begin{equation*}
\begin{bmatrix}
    \mathbf z_{k+1}\\
    \mathbf u_{k+1}
\end{bmatrix} =
    \begin{cases}
        \begin{bmatrix}  % case 1
            \rho\mathbf D\mathbf R\mathbf D^T & \mathbf I-\rho\mathbf D\mathbf R\mathbf D^T\\
            0                   & 0
        \end{bmatrix}
        
        \begin{bmatrix}
            \mathbf z_k\\
            \mathbf u_k
        \end{bmatrix}
        +
        \begin{bmatrix}
            \mathbf D\mathbf R\mathbf y-\frac{\alpha}{\rho}\mathbf 1\\
            \frac{\alpha}{\rho}\mathbf 1
        \end{bmatrix}
        & \text{if } \mathbf D\mathbf x_{k+1}+\mathbf u_k  > \frac{\alpha}{\rho}\mathbf 1,\\
        \\
        \begin{bmatrix}  % case 2
            0                   & 0\\
            \rho\mathbf D\mathbf R\mathbf D^T & \mathbf I-\rho\mathbf D\mathbf R\mathbf D^T 
        \end{bmatrix}
        
        \begin{bmatrix}
            \mathbf z_k\\
            \mathbf u_k
        \end{bmatrix}
        +
        \begin{bmatrix}
            0\\
            \mathbf D\mathbf R\mathbf y
        \end{bmatrix}
        & \text{if } \left| \mathbf D\mathbf x_{k+1}+\mathbf u_k \right| \leq \frac{\alpha}{\rho}\mathbf 1,\\
        \\
        \begin{bmatrix}  % case 3
            \rho\mathbf D\mathbf R\mathbf D^T  & \mathbf I-\rho\mathbf D\mathbf R\mathbf D^T \\
            0 & 0
        \end{bmatrix}
        
        \begin{bmatrix}
            \mathbf z_k\\
            \mathbf u_k
        \end{bmatrix}
        +
        \begin{bmatrix}
            \mathbf D\mathbf R\mathbf y+\frac{\alpha}{\rho}\mathbf 1\\
            -\frac{\alpha}{\rho}\mathbf 1
        \end{bmatrix}
        & \text{if } \mathbf D\mathbf x_{k+1}+\mathbf u_k < -\frac{\alpha}{\rho}\mathbf 1,\\
    \end{cases}
\end{equation*}
where $\mathbf R = (\mathbf I+\rho \mathbf D^T \mathbf D)^{-1}$. Hence, we have
\begin{equation*}
\mathbf q'\left(
\begin{bmatrix}
    \mathbf z_{k+1}\\
    \mathbf u_{k+1}
\end{bmatrix} 
\right)
=
    \begin{cases}
        \begin{bmatrix}  % case 1
            \rho\mathbf D\mathbf R\mathbf D^T & \mathbf I-\rho\mathbf D\mathbf R\mathbf D^T\\
            0                   & 0
        \end{bmatrix}
        & \text{if } |\mathbf D\mathbf x_{k+1}+\mathbf u_k|  > \frac{\alpha}{\rho}\mathbf 1,\\
        \\
        \begin{bmatrix}  % case 2
            0                   & 0\\
            \rho\mathbf D\mathbf R\mathbf D^T & \mathbf I-\rho\mathbf D\mathbf R\mathbf D^T
        \end{bmatrix}
        & \text{if } \left| \mathbf D\mathbf x_{k+1}+\mathbf u_k\right| \leq \frac{\alpha}{\rho}\mathbf 1.\\
    \end{cases}
\end{equation*}
When the conditions on $\mathbf D\mathbf x_{k+1}$ and $\mathbf u_k$ do not fall completely into one category, we have to interpret the above expressions component-wise like in our analysis for the scalar example.

Similarly, for the lasso problem of Section \ref{sec:lasso}
we get \begin{equation*}
\begin{bmatrix}
    \mathbf z_{k+1}\\
    \mathbf u_{k+1}
\end{bmatrix} =
    \begin{cases}
        \begin{bmatrix}  % case 1
            \rho\mathbf R & \mathbf I-\rho\mathbf R\\
            0      & 0
        \end{bmatrix}
        
        \begin{bmatrix}
            \mathbf z_k\\
            \mathbf u_k
        \end{bmatrix}
        +
        \begin{bmatrix}
           \mathbf R\mathbf A^T\mathbf b -\frac{\lambda}{\rho}\mathbf 1\\
            \frac{\lambda}{\rho}\mathbf 1
        \end{bmatrix}
        & \text{if } \mathbf x_{k+1}+\mathbf u_k  > \frac{\alpha}{\rho}\mathbf 1,\\
        \\
        \begin{bmatrix}  % case 2
            0                   & 0\\
            \rho\mathbf R & \mathbf I-\rho\mathbf R
        \end{bmatrix}
        
        \begin{bmatrix}
            \mathbf z_k\\
            \mathbf u_k
        \end{bmatrix}
        +
        \begin{bmatrix}
           0\\
           \mathbf R\mathbf A^T\mathbf b
        \end{bmatrix}
        & \text{if } \left| \mathbf x_{k+1}+\mathbf u_k \right| \leq \frac{\alpha}{\rho}\mathbf 1,\\
        \\
        \begin{bmatrix}  % case 3
            \rho\mathbf R & \mathbf I-\rho\mathbf R\\
            0 & 0
        \end{bmatrix}
        
        \begin{bmatrix}
            \mathbf z_k\\
            \mathbf u_k
        \end{bmatrix}
        +
        \begin{bmatrix}
            \mathbf R\mathbf A^T\mathbf b+\frac{\lambda}{\rho}\mathbf 1\\
            -\frac{\lambda}{\rho}\mathbf 1
        \end{bmatrix}
        & \text{if } \mathbf x_{k+1}+\mathbf u_k < -\frac{\alpha}{\rho}\mathbf 1,\\
    \end{cases}
\end{equation*}
where $\mathbf R = (\mathbf A^T\mathbf A+\rho \mathbf I)^{-1}$, and
\begin{equation*}
\mathbf q'\left(
\begin{bmatrix}
    \mathbf z_{k+1}\\
    \mathbf u_{k+1}
\end{bmatrix} 
\right)
=
    \begin{cases}
        \begin{bmatrix}  % case 1
            \rho\mathbf R & \mathbf I-\rho\mathbf R\\
            0                   & 0
        \end{bmatrix}
        & \text{if } |\mathbf x_{k+1}+\mathbf u_k|  > \frac{\lambda}{\rho}\mathbf 1,\\
        \\
        \begin{bmatrix}  % case 2
            0                   & 0\\
            \rho\mathbf R & \mathbf I-\rho\mathbf R
        \end{bmatrix}
        & \text{if } \left| \mathbf x_{k+1}+\mathbf u_k\right| \leq \frac{\lambda}{\rho}\mathbf 1.\\
    \end{cases}
\end{equation*}

For these multi-dimensional problems with nonsmooth objective function, we find that
the Jacobian of the ADMM iteration function, $\mathbf q'(\mathbf x)$, exists at the solution
$\mathbf x^*$, which is consistent with the observed linear convergence of ADMM
with convergence factor $\rho(\mathbf q'(\mathbf x^*))$.

%\begin{figure}
%\centering
%\caption{Lasso problem (density = 0.06): spectrum of ADMM iteration matrix $\mathbf q'$ and sAA($m$)-ADMM iteration matrix $\mathbf\Psi'$.}
%\label{fig:eigs_M_lasso_3}
%\end{figure}

\end{appendices}
\section*{Declarations}
\textbf{Funding}
HDS gratefully acknowledges support from the Natural Sciences and Engineering Research Council of Canada through the Discovery Grant program (RGPIN-2019-04155).\\
\textbf{Conflicts of interest}
The authors declare that they have no conflict of interest.\\
\textbf{Code availability}
Computer implementation of the algorithms and numerical tests reported on in this paper is freely available at \url{https://github.com/dw-wang/AA-ADMM}.\\
\textbf{Data availability statement}
Data sharing is not applicable to this article as no datasets were generated or analysed.

\end{document}